\documentclass[11pt, reqno]{amsart}
\usepackage[active]{srcltx}
\usepackage {color}
\usepackage{amsmath,amsthm,amssymb}
\usepackage{epsfig}
\usepackage{graphicx}
\usepackage{amssymb}
\usepackage{latexsym}
\usepackage{pdfpages}
\usepackage{bm}
\usepackage{enumitem}
\usepackage{esint}
\usepackage{comment}
\usepackage{float}

\usepackage{ulem} 

\usepackage{ifpdf}
\newcommand{\res}{\!\!\mathop{\hbox{
                                \vrule height 7pt width .5pt depth 0pt
                                \vrule height .5pt width 6pt depth 0pt}}
                                \nolimits}

\def\z{{\bf z}}

\def\divi{\hbox{\rm div\,}}

\ifpdf 
  \usepackage[hidelinks]{hyperref}
\else 
\fi 

 \usepackage[usenames,dvipsnames]{pstricks}
 \usepackage{epsfig}
 \usepackage{pst-grad} 
 \usepackage{pst-plot} 

\usepackage{color}

\usepackage{subcaption}

\newtheorem{theorem}{Theorem}[section]
\newtheorem{lemma}[theorem]{Lemma}
\newtheorem{definition}[theorem]{Definition}
\newtheorem{proposition}[theorem]{Proposition}
\newtheorem{corollary}[theorem]{Corollary}
\newtheorem{remark}[theorem]{Remark}
\newtheorem{example}[theorem]{Example}
\newtheorem*{theorem*}{\it Theorem}

\def\divi{\hbox{\rm div\,}}

\DeclareMathOperator*{\esssup}{ess\,sup}

\def\vint_#1{\mathchoice%
          {\mathop{\kern 0.2em\vrule width 0.6em height 0.69678ex depth -0.58065ex
                  \kern -0.8em \intop}\nolimits_{\kern -0.4em#1}}%
          {\mathop{\kern 0.1em\vrule width 0.5em height 0.69678ex depth -0.60387ex
                  \kern -0.6em \intop}\nolimits_{#1}}%
          {\mathop{\kern 0.1em\vrule width 0.5em height 0.69678ex
              depth -0.60387ex
                  \kern -0.6em \intop}\nolimits_{#1}}%
          {\mathop{\kern 0.1em\vrule width 0.5em height 0.69678ex depth -0.60387ex
                  \kern -0.6em \intop}\nolimits_{#1}}}
\def\vintslides_#1{\mathchoice%
          {\mathop{\kern 0.1em\vrule width 0.5em height 0.697ex depth -0.581ex
                  \kern -0.6em \intop}\nolimits_{\kern -0.4em#1}}%
          {\mathop{\kern 0.1em\vrule width 0.3em height 0.697ex depth -0.604ex
                  \kern -0.4em \intop}\nolimits_{#1}}%
          {\mathop{\kern 0.1em\vrule width 0.3em height 0.697ex depth -0.604ex
                  \kern -0.4em \intop}\nolimits_{#1}}%
          {\mathop{\kern 0.1em\vrule width 0.3em height 0.697ex depth -0.604ex
                  \kern -0.4em \intop}\nolimits_{#1}}}

\def\R{\mathbb R}
\def\N{\mathbb N}
\def\Z{\mathbb Z}
\def\g{\hbox{\bf g}}

\numberwithin{equation}{section}

\def\avint{\mathop{\mathchoice{\,\rlap{-}\!\!\int}
		{\rlap{\raise.15em{\scriptstyle -}}\kern-.2em\int}
		{\rlap{\raise.09em{\scriptscriptstyle -}}\!\int}
		{\rlap{-}\!\int}}\nolimits}

\def\dfrac{\displaystyle\frac}

\def\1{\raisebox{2pt}{\rm{$\chi$}}}

\definecolor{violet(ryb)}{rgb}{0.53, 0.0, 0.69}
\usepackage{mathtools}
\mathtoolsset{showonlyrefs}

 \newcounter{mgncount}

\begin{document}

\title[Anisotropic inverse mean curvature flow]{\bf Weak solutions of Anisotropic (and crystalline) inverse mean curvature flow \\ as limits of $p$-capacitary potentials}

\author[E. Cabezas-Rivas, S. Moll, M. Solera]{  Esther Cabezas-Rivas, Salvador Moll and Marcos Solera}

\address{E. Cabezas-Rivas: Departament de Matem\`atiques,
Universitat de Val\`encia, Dr. Moliner 50, 46100 Burjassot, Spain.
 {\tt esther.cabezas-rivas@uv.es  }}

\address{S. Moll: Departament d'An\`{a}lisi Matem\`atica,
Universitat de Val\`encia, Dr. Moliner 50, 46100 Burjassot, Spain.
 {\tt j.salvador.moll@uv.es }}
\address{M. Solera: Departament d'An\`{a}lisi Matem\`atica,
Universitat de Val\`encia, Dr. Moliner 50, 46100 Burjassot, Spain. {\tt marcos.solera@uv.es } \vspace{0.1cm}\newline
\phantom{MMM. Solera:} Departamento de Matem\'aticas, Universidad Aut\'onoma de Madrid, C/
Francisco Tom\'as y Valiente, 7, Facultad de Ciencias, m\'odulo 17, 28049 Madrid, Spain.
}


\date{\today. This work is partially supported by project PID2022-136589NB-I00 and the network RED2022-134077  funded by MCIN/AEI/10.13039/501100011033 and by ERDF  A way of making Europe. The  first author is partially supported by project PID2019-105019GB-C21 funded by MCIN/AEI/10.13039/501100011033/ and by ERDF  A way of making Europe; and by Conselleria d'Educaci\'o, Universitats i Ocupaci\'o, project AICO/2021/252. The second and third authors have been partially supported  by Conselleria d'Educaci\'o, Universitats i Ocupaci\'o, project AICO/2021/223, as well as the network RED2022-134784-T funded by MCIN/AEI/10.13039/501100011033. The third author has also been supported by the Conselleria d'Educaci\'o, Universitats i Ocupaci\'o, programme ``Subvenciones para la contratación de personal investigador en fase postdoctoral'' (APOSTD 2022),  Ref. CIAPOS/2021/28; the MICINN (Spain) and ERDF, project PID2021-124195NB-C32; and by the ERC Advanced Grant 834728.}


\setcounter{tocdepth}{1}

%

\begin{abstract}
We construct weak solutions of the anisotropic inverse mean curvature flow (A-IMCF) under very mild assumptions both on the anisotropy (which is simply a norm in $\R^N$ with no ellip\-ticity nor smoothness requirements, in order to include the crystalline case) and on the initial data. By means of an approximation procedure introduced by Moser, our solutions are limits of anisotropic $p$-harmonic functions or $p$-capacitary functions (after a change of variable), and we get uniqueness both for the approximating solutions (i.e., uniqueness of $p$-capacitary functions) and the limiting ones. Our notion of weak solution still recovers variational and geometric definitions similar to those introduced by Huisken-Ilmanen, but requires to work within the broader setting of $BV$-functions. Despite of this,  we still reach classical results like the continuity and exponential growth of perimeter, as well as outward minimizing properties of the sublevel sets. Moreover, by assuming the extra regularity given by an interior rolling ball condition (where a sliding Wulff shape plays the role of a ball), the solutions are shown to be continuous and satisfy Harnack inequalities. Finally, examples of explicit solutions are built.

\end{abstract}

\maketitle



{ \renewcommand\contentsname{Contents }
\setcounter{tocdepth}{3}
 }

\vspace*{-1cm}

 \section{Introduction and statement of main results}

 \subsection{Moser's approach to inverse mean curvature flow (IMCF)}

For $p >1$, let $u$ be a positive $p$-harmonic function or $p$-capacitary potential, that is, a solution to
\begin{equation} \label{p-lap}
{\rm div}\big(|\nabla u|^{p-2} \nabla u \big) = 0.
\end{equation}
After the change of variable
\begin{equation} \label{change-var}
v = (1-p) \log u,
\end{equation}
the new function $v$ satisfies
\begin{equation} \label{p-IMCF}
{\rm div}\big(|\nabla v|^{p-2} \nabla v \big) = |\nabla v|^p.
\end{equation}

 R. Moser developed in \cite{Moser, Moser08} the striking idea of taking limits as $p\searrow 1$  of solutions $v_p$ to \eqref{p-IMCF} leading to  weak solutions of the degenerate elliptic equation:
 \begin{equation} \label{level-IMCF}
 {\rm div}\left(\frac{\nabla v}{|\nabla v|} \right) = |\nabla v|,
 \end{equation}
which has its own independent interest since, if $v$ is smooth and $\nabla v \neq 0$, then the level sets of $v$ evolve according to the IMCF, that is, each point moves in the unit outward normal direction with velocity given by  the inverse $1/H$ of the mean curvature.

In particular, Moser constructs weak solutions of \eqref{level-IMCF}  by using
 explicit barriers to first solve \eqref{p-lap} with suitable Dirichlet boundary conditions.  Unlike \eqref{level-IMCF}, the $p$-Laplace equation \eqref{p-lap} has the advantage of being the Euler-Lagrange equation of the variational problem of minimizing the $p$-capacity.

Roughly speaking, the overall goal of this paper is to reproduce Moser's successful technique within an anisotropic framework by imposing the mildest possible requirements both for the anisotropy (in order to include the crystalline case) and the domains where the solutions are defined.  Moving beyond traditional Euclidean isotropic settings is essential to allow for situations where attributes vary with direction, as this additional flexibility will be decisive for applications to crystal development or noise reduction.

Having this in mind, as a starting point, in \cite{CRMS} we studied existence and regularity properties of solutions to the anisotropic version of \eqref{p-lap}.

\subsection{Anisotropic implementation of Moser's technique}

 Let $N \geq 2$ and $\Omega \subset \R^N$ be an open bounded set with Lipschitz-continuous boundary. Hereafter, we will use the notation $\Omega^e$ for the exterior domain $\R^N \setminus \overline{\Omega}$. We can formally define the boundary value problem:
 \begin{equation}\label{EquationU}
 \left\{ \begin{array}{ll} {\rm div}\big( F^{p-1}(\nabla u) \partial F(\nabla u)\big) = 0 & \text{in }  \Omega^e, \smallskip \\
 u = 1 & \text{on }  \partial\Omega, \smallskip \\
 u = 0 & \text{as }  |x| \to \infty, \end{array} \right.
 \end{equation}
 where $F$ is a norm on $\R^N$. Here the
 equality on the boundary is understood in the sense of traces, and $\partial F$ stands for an element of the subdifferential (for the sake of clarity, we keep this formal notation to present our main results, see subsections \ref{def-sol-2} and \ref{def-sol1} for precise definitions).

 As before, the transformation \eqref{change-var}, leads to
 \begin{equation}\label{EquationV}
 \left\{ \begin{array}{ll} {\rm div}\big( F^{p-1}(\nabla v) \partial F(\nabla v)\big) = F^p(\nabla v) & \text{in }  \Omega^e, \smallskip \\
 v = 0 & \text{on }  \partial\Omega, \smallskip \\
 v \to  \infty & \text{as }  |x| \to \infty. \end{array} \right.
 \end{equation}

\noindent As Moser's scheme suggests,  the natural candidate for a definition of anisotropic or crystalline IMCF (A-IMCF) should be the limit of \eqref{EquationV} as $p \searrow 1$:
 \begin{equation}\label{EquationVp1}
 \left\{ \begin{array}{ll} {\rm div}\big(\partial F(D v)\big) = F(D v) & \text{in }  \Omega^e, \smallskip \\
 v = 0 & \text{on }  \partial\Omega, \smallskip \\
 v \to  \infty & \text{as }  |x| \to \infty, \end{array} \right.
 \end{equation}
 where $D v$ represents the distributional derivative of a $BV$-function, which makes sense as a Radon measure (for a comprehensive introduction to this framework, see subsection \ref{BV-section}).

 Indeed, our main result confirms the above heuristics.
 	\begin{theorem}[Existence \& uniqueness of weak solutions of A-IMCF] \label{umin-trapped}
 	Let $1<p<N$ and consider  $\Omega \subset\R^N$ a bounded domain with Lipschitz  boundary. If $u_p$ solves \eqref{EquationU} and we set $v_p = (1-p) \log u_p$,   there exists a sequence $p_k \searrow 1$ and a function  $v  \in BV_{\rm loc}(\Omega^e)\cap L^\infty_{\rm loc}(\Omega^e)$ such that $v_{p_k} \to v$ in $L^q_{\rm loc}(\Omega^e)$ for every $1 \leq q <\frac{N}{N-1}$ and the limit has the following properties:
 	\begin{enumerate}
 		\item[{\rm (a)}] $v$ is the unique weak solution of \eqref{EquationVp1}.
 		\item[{\rm (b)}] The discontinuity or singular set of $v$ has measure zero with respect to the $(N-1)$-dimensional Hausdorff measure.
 		\item[{\rm (c)}] We can find constants $0 < r_1 < r_2$ with $\mathcal W_{r_1} \subset \Omega \subset  \mathcal W_{r_2}$ so that
 		\[v_{r_1} \leq v \leq v_{r_2}, \quad \text{where } \ v_{r_i} = (N-1) \log\bigg(F^\circ\Big(\frac{\cdot}{r_i}\Big)\bigg), \quad i = 1,2.\]
 		Here $F^\circ$ denotes the dual of the norm $F$ and $\mathcal W_r$ its corresponding ball (Wulff shape) of radius $r$.
 	\end{enumerate}
 \end{theorem}
\noindent Observe that, by translation invariance, we can suppose that $0\in\Omega$ and, since $\Omega$ is bounded, there exist $0<r_1<r_2$ such that $\mathcal W_{r_1}\subset\Omega\subset\mathcal W_{r_2}$.

 Up to our knowledge, the most general existence and uniqueness result in this spirit can be found in \cite{IAMCF}. The authors obtain a weak solution  $v\in  \rm Lip_{\rm loc}(\Omega^e)$ in the case that $\Omega$ has smooth boundary for norms $F\in C^{\infty}(\R^N\setminus\{0\})$ with a strictly positive definite Hessian matrix.

 Here, we broaden the existing literature's findings in a number of ways. In fact, we just require the anisotropy to be a norm, without making any additional commitments about smoothness or uniform ellipticity. Specifically, we permit norms whose dual unit balls have corners and/or flat pieces, including crystalline situations as the $\ell_\infty$ or $\ell_1$ norms in $\R^N$. Moreover, about the regularity of $\Omega$, we just ask for Lipschitz continuity of the boundary, instead of the traditional $C^\infty$ much stronger constraint.

 In short, we deal with an unfriendly setting because we will not be able to produce any argument that involves gradient or second derivatives of the functions $v$ nor principal curvatures (even in weak sense) of $\partial \Omega$. To overcome these additional technical difficulties, the proof of the existence in \autoref{umin-trapped} strongly relies on the ideas developed in \cite{ABCM} to study a Dirichlet problem coming from the {\it total variation flow}, and the subsequent adaptation of them to study IMCF with a damping term on bounded \cite{MazSer13} and unbounded \cite{MazSer15} Euclidean domains.

 More precisely, we start by introducing a notion of weak solution to \eqref{EquationVp1} (see Definition \ref{weak-IMCF-def}), where the equation makes sense as an equality for measures and which relies on the existence of a bounded vector field $\z$, whose divergence is a Radon measure that plays the role of the left hand side of \eqref{level-IMCF}, even if $\nabla v$ may vanish. To overcome this, we exploit intensively Anzellotti’s theory of pairings \cite{Anz} and the total variation adapted to anisotropies by \cite{AmBe}. Let us highlight that measure divergence fields have already shown their key role in physics as models of capillarity for perfectly wetting fluids \cite{LeSa} and in continuum mechanics \cite{Stu}.

 While we managed to follow the existence proof of \cite{MazSer13, MazSer15} with suitable adaptations and correction of several issues, the comparison principles there leading to uniqueness fail in the case that $F$ is not strictly convex. Luckily, we manage to combine a trick inspired by the uniqueness proof by Huisken and Ilmanen \cite[Theorem 2.2]{HuIl} with the use of truncation functions to achieve comparison results (and hence uniqueness)  for the A-IMCF \eqref{EquationVp1}. We stress that our uniqueness claim and proof differs quite from that in \cite{HuIl} because of our different notion of weak solution.

 Furthermore, these new insights also work for the approximating problems \eqref{EquationU} and \eqref{EquationV}, even for generic boundary conditions, meaning that as a by-product of the techniques developed here we attain an upgrade of our existence and regularity results in \cite{CRMS}. With more details,
 \begin{corollary} \label{uniq-CRMS1}
 	Let  $1 < p < \infty$, $\Omega \subset \R^N$ a bounded domain with Lipschitz boundary and $\varphi \in W^{1-\frac1{p},p}(\partial \Omega)$. Consider either $D = \Omega$ or $D = \Omega^e$, then the weak solution of
 	\begin{equation}\label{EquationUvarphi}
 	\left\{ \begin{array}{ll} {\rm div}\big( F^{p-1}(\nabla u) \partial F(\nabla u)\big) = 0 & \text{in }  D\smallskip \\
 	u = \varphi & \text{on }  \partial \Omega, \smallskip \\
 	u \to 0 \text{ as } |x| \to \infty & \text{if }  D = \Omega^e, \quad p < N\end{array} \right.
 	\end{equation}
 	constructed in \cite{CRMS} is unique. In particular, if $\varphi \equiv 1$, then uniqueness of minimizers  in $\{u\in L^{pN/(N-p)}(\R^N)\,:\, \nabla u \in L^p(\R^N), \, u\geq 1 \text{ in } \Omega\}$ of the {\it anisotropic p-capacity functional}
 	\begin{equation} \label{pcap-def}
 	{\rm Cap}^F_p(\overline \Omega):=\inf\left\{\int_{\R^N} F^p(\nabla u)\,\mathrm{d}x : u\in C_c^\infty(\R^N), u\geq 1 {\rm \ in \ } \overline{\Omega}\right\}.
 	\end{equation}
 	 is guaranteed.
 \end{corollary}
\noindent Recall that physically speaking those minimizers are known as $p$-capacitary potentials and are related to thermal or electric conductivity problems in heterogeneous media, where properties of the ambient are allowed to depend upon the direction.

\subsection{Further regularity and Harnack inequalities}
	 In \cite{CRMS} we proved that solutions to \eqref{EquationU} are Lipschitz continuous provided that the domain is regular enough in the following sense:
	\begin{definition}[Uniform interior ball condition] \label{UIBC}
	Let $r>0$. We say that $\Omega$ satisfies the $\mathcal W_r$-condition if, for any $x\in\partial \Omega$, there exists $y\in \R^N$ such that
	\[\mathcal W_r+y\subseteq \overline{\Omega} \quad \text{ and } \quad x\in\partial \left(\mathcal W_r+y\right).\]
\end{definition}
Accordingly, it is natural to expect also extra regularity for the solutions of \eqref{EquationVp1} under the same circumstances. Indeed,
\begin{theorem} \label{cont+Harnack:thm}
	Let $r>0$ and suppose that $\Omega$ satisfies the $\mathcal W_r$-condition. Then there exists  $v  \in BV_{\rm loc}(\Omega^e)\cap C^0(\Omega^e)$ unique solution of \eqref{EquationVp1}. Moreover, one can find a constant $\mathcal C>0$ depending only on the dimension $N$, such that, if $x_0\in\Omega^e$ and $\rho>0$ satisfy $\mathcal W_{2 \rho}(x_0)\subset \Omega^e$, then
	\[{\rm osc}_{\mathcal W_\rho(x_0)} v \leq \mathcal C,\]
	 where ${\rm osc}_B v := \sup_B v - \inf_B v$ is the oscillation of $v$ on a subset $B$ .
\end{theorem}
\noindent The proof follows the strategy in \cite{Moser08} of constructing explicit barrier functions leading to continuity via comparison arguments and getting a Harnack inequality of the type
\[\sup_{\mathcal W_\rho(x_0)} u_p \leq e^{\frac{\mathcal C}{p-1}} \inf_{\mathcal W_\rho(x_0)} u_p\]
for solutions to \eqref{EquationU}, which after the typical transformation \eqref{change-var} and passing to the limit as $p\searrow 1$, yields the claimed estimate for $v$.

\subsection{\lq\lq Classical" geometric properties of solutions}
On the other hand, despite of our milder regularity assumptions and our a priori different notion of solutions, the output of the developments here does not differ as much as could be expected from those by Huisken-Ilmanen, since our solution $v$ has also a similar variational nature (see Proposition \ref{def-var}) and its corresponding sublevel sets
\begin{equation} \label{def-sublevel}
E_t = \{x \,: \, v(x) < t\} \quad \text{and } \quad G_t = \{x \,: \, v(x) \leq t\}
\end{equation}
also minimize a suitable functional (cf. Proposition \ref{equiv-sol}). Additionally, we recover the continuity of the anisotropic version of the perimeter, exponential growth and minimizing hull properties stated in \cite{HuIl}.
\begin{theorem}[Geometric properties of solutions of A-IMCF] \label{sublevels-all} The sublevel sets in \eqref{def-sublevel} satisfy the following properties:
\begin{enumerate}
	\item[{\rm (a)}] $G_\tau$ and $E_t$ are sets of finite perimeter for every $t > \tau \geq 0$.
	
	\item[{\rm (b)}] The anisotropic perimeter $P_{_F}$ of $G_t$ and $E_t$ coincide for all $t >0$.
	
	\item[{\rm (c)}] $e^{-t}  P_{_F}(E_t)$ is constant for $t > 0$.
	
	\item[{\rm (d)}] $E_t$ is outward $F$-minimizing for any $t > 0$, that is, it  minimizes the anisotropic perimeter when compared with all bounded subsets that contain it, while enclosing $\Omega$.

 \item[{\rm (e)}] $G_t$ is strictly outward $F$-minimizing for any $t\geq 0$. Moreover, among the anisotropic perimeter minimising envelopes of $\Omega$, $G_0$ is the one that maximises the Lebesgue measure.
\end{enumerate}	
\end{theorem}
We actually recover these {\it classical} properties (see  \cite[Property 1.4 and Lemma 1.6]{HuIl}) under much weaker regularity assumptions. This makes the proof quite more convoluted since our functions are not continuous and, therefore, claims that were immediate in \cite{HuIl}, require a more careful study here, including the use of notions like  the measure theoretic interior or reduced boundary, not present in the original developments.

For the proof of (a) we follow \cite[Section 4]{MazSer15} and fill some gaps of the arguments there, while adapting them to the anisotropic environment. For the rest of items, we argue in the spirit of \cite{HuIl} but modifying appropriately some notions to better suit our purposes (see section \ref{out-min} for precise definitions and further geometric properties of solutions).

\subsection{Relevance of IMCF}   Given a smooth manifold $M$ of dimension $N-1$, a classical solution of the IMCF is a family of hypersurfaces $X: M\times [0, T) \rightarrow  \R^N$ satisfying
\begin{equation} \label{IMCF}
\partial_t X(p,t) = \frac{\nu}{H}(p,t), \qquad p \in M, \quad  0<t<T,
\end{equation}
where $\nu$ is the unit outward normal to the evolving hypersurfaces $M_t = X(M, t)$ and $H = {\rm div} \, \nu$ is the corresponding mean curvature, which is assu\-med to be positive everywhere so that the above equation is well-defined.

If the initial condition is smooth, star-shaped and mean convex (with $H > 0$), then Gerhardt \cite{Gerhardt} and Urbas \cite{Urbas} obtained independently that the solution exists for all times and the rescaled hypersurfaces $e^{-\frac{t}{N-1}} M_t$ converge to a round sphere as $t \to \infty$. The anisotropic version of this result for smooth anisotropy and initial condition was proved in \cite{Xia2}.

However, without geometric conditions on the initial data, singularities may develop; in fact, for the case $N =3$ it is proved in \cite{Smo} that singularities happen only if $\inf_{M_t} H \to 0$ during the flow. Because of this, Huisken and Ilmanen \cite{HuIl} introduced the notion of weak solution to \eqref{IMCF} using the aforementioned level-set approach. Later on, they showed higher regularity properties in \cite{HuIl2}: every weak solution to the flow is smooth after the first instant $t_0$ where a level set $\partial E_{t_0}$ becomes star-shaped.

The reason behind the fact that this flow attracts so much attention is its well established effectiveness  to generate interesting inequalities which are geometric and physically meaningful. In a nutshell, it was applied by Huisken-Ilmanen \cite{HuIl} to prove Riemannian Penrose inequality in asymptotically flat 3-manifolds, from which one gets a lower bound for the total energy of a universe
 at a space-like infinity as a function of the largest black hole in the system. Later on, Brendle–Hung–Wang \cite{BHW}
derived a Minkowski type inequality for Anti-de Sitter Schwarzschild manifolds by means of IMCF, which leads to a Gibbons–Penrose’s inequality in Schwarzschild spacetime \cite{BW}. Further development of such a a powerful tool in anisotropic situation should lead to new useful inequalities (see e.g. anisotropic Minkowski inequalities in \cite{Xia2}).

\subsection{Structure of the paper} The paper is organized as follows. We first
introduce in section \ref{back} the basic background material about anisotropies, functions of bounded variation, divergence-measure vector fields, pairing of measures and anisotropic total variation; in particular, the latter requires to include the technical details of the proofs that lead us to get Cauchy-Schwarz type inequalities (Lemma \ref{prop-ATV-2}) under milder assumptions than those in the previous literature. In section \ref{sec-def-flow} we give our notion of solution of A-IMCF (see Definition \ref{weak-IMCF-def}) and deduce the first consequences, that will be crucial for the rest of the paper, like the nullity of the singular set with respect of the Hausdorff measure (Lemma \ref{no-jump}). Next, we prove a comparison result (Theorem \ref{unique}) for sub- and supersolutions of A-IMCF, leading to uniqueness of solutions. We start with Moser's strategy in section \ref{Moser-aprox}, where we define the corresponding approximating problems (before and after the change of variable) and obtain suitable comparison theorems for them; this leads us to attain Corollary \ref{uniq-CRMS1}. In turn, the proofs of Theorem \ref{umin-trapped} and Theorem \ref{cont+Harnack:thm} are carried out in sections  \ref{existence-sect} and \ref{cont-harnack}, respectively. The next natural step is to worry about the geometric properties of the sublevel sets of the solutions of A-IMCF, which stands for the content of sections \ref{sublevel-1} and \ref{out-min}; these altogether provide a proof of Theorem \ref{sublevels-all}. Finally, we work out explicit examples of solutions: the expected evolution of Wulff shapes, a rectangle whose vertices do not growth linearly, as well as a concrete case where {\it fattening} occurs.

 \section{Background material: anisotropies and BV functions} \label{back}

 \subsection{Anisotropies}

A  function $F:\R^N \rightarrow [0,+\infty[$ is said to be an {\it anisotropy} if it is convex, positively 1-homogeneous (i.e., $F(\lambda x) = \lambda F(x)$ for all $\lambda >0$ and all $x \in \R^N$) and coercive. We will always consider additionally that $F$ is even, that is, a norm. In particular, as all norms in $\R^N$ are equivalent, 
there exist  constants $0<c\le C<\infty$ such that
\begin{equation}\label{FequivEuclid}
c\vert \xi\vert \le F(\xi) \le C \vert\xi\vert
\end{equation}
where $\vert \cdot \vert$ is the Euclidean norm in $\R^N$.

We define the dual or polar function $F^\circ:\R^N\rightarrow [0,+\infty[$ of $F$ by
\[F^\circ(\xi):=\sup \{ \xi\cdot\xi^\star \, : \, \xi^\star\in\R^N, \, F(\xi^\star)\le 1\}
=\sup \left\{\frac{\xi \cdot \xi^\star}{F(\xi^\star)}\, : \, \xi^\star\in\R^N \right\}\]
for every $\xi\in\R^N.$
It can be verified that $F^\circ$ is convex, lower semi-continuous and $1$-positively homogeneous. Moreover, \eqref{FequivEuclid} leads to
\begin{equation}\label{F0equivEuclid}
\frac1{C}\vert \xi\vert \le F^\circ(\xi) \le \frac1{c} \vert\xi\vert.
\end{equation}

From the definition of $F^\circ$ one gets a Cauchy-Schwarz type inequality of the form
\begin{equation} \label{CS-ineq}
x \cdot \xi \leq F^\circ(x) F(\xi) \qquad \text{ for all } \ x, \xi \in \R^N.
\end{equation}

The Wulff shape $\mathcal W_r$ of $F$ (centred at $0$) is defined by
$$\mathcal W_r:=\{ \xi^\star\in\R^N \, : \, F^\circ(\xi^\star)<r\}.$$
As we are dealing with even anisotropies, $\mathcal W_r$ is a centrally symmetric convex body. We say that $F$ is {\it crystalline} if, furthermore, $	\mathcal W_r$ is a convex polytope.  More generally, the Wulff shape centred at $x_0$ of radius $r$ is given by
\[x_0 + \mathcal W_r = \mathcal W_r(x_0) = \{x \in \R^N \, : \, F^\circ(x-x_0) = r\}.\]

 Since both $F$ and $F^\circ$ are norms, it follows that they are Lipschitz functions. Therefore, there exists $\nabla F$ and $\nabla F^\circ$ a.e.~in $\R^N$.  Hence, arguing on points where $\nabla F^\circ(x)$ is well-defined, the following generalization of \cite[Proposition 5.3]{Xia} holds (with the same proof):
\begin{lemma}\label{relationsF} Let $F$ be any norm in $\R^N$, then for a.e. $x \in \R^N \setminus\{0\}$,
\begin{itemize}
  \item[{\rm ($i$)}] $F(\nabla F^\circ (x))=1.$
  \item[{\rm ($ii$)}] $\displaystyle \frac{x}{F^\circ (x)}\in \partial F(\nabla F^\circ(x))$.
\end{itemize}
  \end{lemma}

\subsection{Approximate continuity and differentiability of $L^1$ functions}

Throughout the paper, set $N \geq 2$ a fixed integer, $\mathcal L^N$ and $\mathcal H^{N-1}$ will denote the Lebesgue measure and $(N-1)$-dimensional Hausdorff measure in $\R^N$, respectively. $|E|$ means the $\mathcal L^N$-Lebesgue measure of a subset $E \subset \R^N$ and $\fint_E = \frac1{|E|} \int_E$ the average integral.  Moreover, the integrals on the boundary of a set are computed with respect to $d \mathcal H^{N-1}$.

Hereafter, $U$ stands for any open subset of $\R^N$ with Lipschitz continuous boundary and, as usual,  $B_r(x)$ represents the Euclidean ball in $\R^N$ centred at $x$  with radius $r$. Additionally, we consider the two half balls determined by a given direction $\nu$:
\[B_r(x, \nu)^{\pm} :=\big\{y \in B_r(x) \, : \, \pm (y -x) \cdot \nu >0\big\},\]
where $\cdot$ denotes the Euclidean inner product, and we will use $|u|$ for the Euclidean norm of a function $u$.

\begin{definition}[Approximate discontinuity and jump sets] \label{lebesgue} Let $u \in L^1_{\rm loc}(U)$.
	
	{\tiny $\bullet$} The approximate discontinuity set $S_u \subset U$ is defined in terms of its complement, i.e., $x \not \in S_u$ if there exists a unique $\tilde u(x) \in \R$ such that
	\[\lim_{r \searrow 0}  \fint_{B_r(x)} \big|u(y) - \tilde u(x)\big| \, dy = 0.\]
	
	{\tiny $\bullet$} $u$
	is approximately continuous at $x$ if $x \not\in S_u$ and $\tilde u(x) = u(x)$, that is, if $x$ is a Lebesgue point of $u$.
	
	{\tiny $\bullet$} We now distinguish jump points, where the nature of the discontinuity is more specific: $x \in J_u$ is an approximate jump point of $u$ if there exist $u^+(x) \neq u^-(x) \in \R$ and $\nu_u(x) \in \mathbb S^{N-1}$ such that
	\[\lim_{r \searrow 0}  \fint_{B_r(x, \nu_u(x))^{\pm}} \big|u(y) - u^\pm(x)\big| \, dy = 0.\]
	
\end{definition}

By Lebesgue's differentiation theorem, one can show (cf.~\cite[Proposition 3.64 (a)]{Ambrosio}) that $\tilde u$
exists a.e.~in $U$, and coincides with $u$; in particular, $|S_u| = 0$.

Now take a family of mollifiers $\big\{\rho_\varepsilon(x)= \frac1{\varepsilon^N} \rho(x/\varepsilon)\big\}_{\varepsilon >0}$, where $\rho \in C^\infty_c(\R^N)$ is a non-negative, even function, with ${\rm supp} \rho \subset B_1(0)$ and $\int_{\R^N} \rho = 1$. By \cite[Propositions 3.64 (b) and 3.69 (b)]{Ambrosio}, we have the following result:

\begin{lemma}[Convergence to the precise representative] \label{conv-u-ast} For any $u \in L^1_{\rm loc}(U)$, the mollified functions $u \ast \rho_\varepsilon \in C^\infty(U)$ pointwise converge as $\varepsilon \searrow 0$ out of $S_u \setminus J_u$ to the precise representative $u^\ast: U \setminus(S_u \setminus J_u) \rightarrow \R$ given by
	\begin{equation} \label{def-uast}
	u^{\ast}(x):=\left\{\begin{array}{ll} \tilde u(x), & \text{if } x \in U \setminus S_u \medskip \\ \dfrac{u^+(x)+u^-(x)}{2}, & \text{if } x \in J_u. \end{array} \right.
	\end{equation}
\end{lemma}

There is also a notion of approximate differentiability.
\begin{definition}
$u \in L^1_{\rm loc}(U)$ is said to be approximate differentiable at $x \in U \setminus S_u$ if there exists $\nabla u(x) \in \R^N$, called approximate gradient of $u$ at $x$, such that
\[\lim_{r \searrow 0} \frac1{r} \fint_{B_r(x)} \big|u(y) - \tilde u(x) - \nabla u(x) \cdot (y-x)\big| \, dy = 0.\]
\end{definition}

\subsection{Subdifferential of a functional and $F^p$}

Let $X$ be a reflexive Banach space with dual $X'$, and denote by $\big< \,\cdot\,, \,\cdot\,\big>$ the
pairing between $X'$ and $X$.
A $\mathcal A: X \rightarrow 2^{X'}$  multivalued operator on $X$
	$\mathcal A$ is said to be {\it monotone} if
	\[\big<\tilde \xi - \tilde\eta, \xi - \eta\big> \geq 0 \qquad \text{for every} \quad  (\xi, \tilde\xi), (\eta, \tilde\eta) \in {\rm graph}(\mathcal A).\]

Let  $\Phi: X \rightarrow \R \cup \{\infty\}$ be  lower semicontinuous, proper and convex; its subdifferential $\partial \Phi: X \rightarrow 2^{X'}$ is a multivalued operator given as follows:
\begin{equation} \label{orig-def-subD}
\delta \in \partial \Phi(\zeta) \quad \Longleftrightarrow \quad \Phi(\eta) - \Phi(\zeta) \geq \big<\delta, \eta-\zeta\big> \quad \text{for all } \eta \in X.
\end{equation}
This is a well-known example of a monotone operator from $X$ to $X'$.

In particular, we are interested in $\partial F^p$, which can be written in terms of $\partial F$ by means of the chain rule for subdifferentials (see \cite[Corollary 16.72]{BC}):
\begin{equation} \label{charFp}
\z \in \partial F^p(\xi) \quad \Longleftrightarrow \quad
	\z = p F^{p-1}(\xi) \tilde \z, \quad \text{with} \quad \tilde \z \in \partial F(\xi).
\end{equation}
This is indeed very useful because there is a well-known characterization of the subdifferential for 1-homogeneous functions:
\begin{equation} \label{char-subdif-1hom}
\tilde \z \in \partial F(\nabla v) \quad \text{ if, and only if,} \quad  \tilde \z \cdot \nabla v = F(\nabla v) \quad \text{ and } \quad F^\circ(\tilde \z) \leq 1.
\end{equation}

\subsection{Functions of bounded variation} \label{BV-section} For a detailed study of the space of functions of bounded variation, we refer to \cite{Ambrosio}.
Let us denote by $\mathfrak{M}(U; \R^N)$ the space of all $\R^N$-valued Radon measures, which is isomorphic to the product space $\mathfrak{M}(U)^N$.
\begin{definition}
We say that $u \in L^1(U)$ is a function of bounded variation, and write $u \in {\rm BV}(U)$, if its distributional derivative $Du$, defined by
\[\big<Du, \varphi\big> := - \int_U u \, {\rm div} \varphi \, dx \quad \text{for all} \quad \varphi \in C^1_c(U; \R^N),\]
belongs to $\mathfrak{M}(U; \R^N)$. Alternatively, one can give the following equivalent characterization: $u \in BV(U)$ if, and only if, $u\in L^1(U)$ and the total variation of $Du$ is finite, meaning that
\begin{equation}
|Du|(U) = \int_U d  |Du| := \sup\bigg\{\int_U u \, {\rm div} \varphi \, : \, \varphi \in C^1_c(U; \R^N), \, \|\varphi\|_\infty \leq 1\bigg\} < \infty.
\end{equation}
\end{definition}
\noindent Unlike Sobolev spaces, $BV(U)$ includes characteristic functions
of sufficiently regular sets. The next result gathers some properties of $BV$-functions that will be useful in the sequel:
\begin{lemma} \label{jump-lem}

{\rm (a)} \cite[Theorem 3.78]{Ambrosio} Every $u\in BV(U)$ is approximately differentiable $\mathcal L^N$-a.e.\ in $U$ and $S_u$ is countably $\mathcal H^{N-1}$-rectifiable; that is,
	\[S_u = (S_u \setminus J_u) \cup \bigcup_{i \in \mathbb N} K_i, \quad \text{where} \quad \mathcal H^{N-1}(S_u \setminus J_u) = 0,\]
and each $K_i$ is a compact subset of a  $(n-1)$-dimensional Lipschitz manifold.

{\rm (b)} $BV(U)$ is a Banach space when it is endowed with the norm
\[\|u\|_{BV(U)}:= \|u\|_{L^1(U)} + |Du|(U).\]

\end{lemma}

\subsection{Decomposition of the measure $Du$} Given $\mu, \tilde \mu$ two positive measures, we say that  $\tilde \mu$ is absolutely continuous with respect to $\mu$, and write $\tilde \mu \ll  \mu$, if, given a measurable set $E$, $\mu(E) = 0$ implies $\tilde \mu(E) = 0$. In this case,  $\frac{\tilde \mu}{\mu}$  means the Radon–Nikodým derivative of $\tilde \mu$ with respect to $\mu$ (called density when $\mu$ is a Lebesgue measure). In turn, $\mu \res E$ is the restriction of $\mu$ to the set $E$.

Notice that Lemma \ref{jump-lem} (a) implies that it is possible to define for  $\mathcal H^{N-1}$-a.e. $x \in J_u$ a unit normal vector field  to $J_u$; in fact, $\nu_u(x)$ is an orientation of $J_u$, and coincides with $\frac{Du}{|Du|}(x)$. Moreover, by Lemma \ref{conv-u-ast}, we get that the mollified functions $u \ast \rho_\varepsilon$, with $u \in BV(U)$, converge to the precise representative $u^\ast$ for $\mathcal H^{N-1}$-a.e.\ point in the domain of $u$.

We may decompose $Du$ in three mutually singular measures (see e.g.\cite[Definition 1.24]{Ambrosio}), as follows (see \cite[Section 3.9]{Ambrosio}):
\begin{equation} \label{Du-dec}
Du = \nabla u \, \mathcal L^N + D^s u, \, \text{with} \, \left\{\begin{array}{l} D^s u = D^j u + D^c u, \medskip \\ D^j u  = (u^+ - u^-)\, \nu_u \, \mathcal H^{N-1}\res J_u, \medskip \\ D^c u = D^s u\res ( U \setminus S_u).\end{array}\right.	
\end{equation}
Here $D^s u$ is the singular part of $Du$  with respect to Lebesgue measure $\mathcal L^N$, and  for $\mathcal L^N$-a.e.\ $x \in U$ the vector $\nabla u(x)$ is the
approximate gradient of $u$ at $x$ (cf.\ \cite[Theorem 3.83]{Ambrosio}). $D^j u$ and $D^c u$ are measures known as the jump and Cantor part of $Du$, respectively. Sometimes, one call  $D^d = \nabla u \, \mathcal L^N  + D^c$ the diffuse part.

\subsection{Divergence-measure fields and weak normal traces}

Here we revisit the integration by parts formula with weak regularity assumptions for the domains of integration and the functions involved. With this aim, consider the space of vector fields whose divergence in the sense of distributions is a bounded Radon measure, i.e.,
\[\mathcal D \mathfrak M^\infty(U)= \{\z \in L^\infty(U; \R^N)\, :\, {\rm div}\, \z \in \mathfrak{M}(U; \R^N) \text{ and is finite}\},\]
where ${\rm div} \, \z$ acts as a measure as follows:
\[\int_U \varphi \, d({\rm div \, \z}) = -\int_U \z \cdot \nabla \varphi \, dx \quad \text{for all} \ \varphi \in C^\infty_c(U).\]

Let $B \subset U$ any bounded open set with Lipschitz continuous boundary and $\nu$ the corresponding outward unit normal on $\partial B$. Now it makes sense to define as in \cite[Theorem 1.2]{Anz} a weak trace of the normal component of $\z \in \mathcal D \mathfrak M^\infty_{\rm loc}(U)$  on $\partial B$ as the linear operator
\[[\,\cdot\,, \nu]: \mathcal D \mathfrak M^\infty(B) \rightarrow L^\infty(\partial B) \quad \text{such that} \quad \big\|[\z, \nu] \big\|_{L^\infty(\partial B)} \leq \|\z\|_\infty,\]
which coincides with $\z (x)\cdot \nu(x)$ for all $x \in \partial B$ when $\z$ is smooth enough.

Here we gather the main properties of divergence-measure fields:
\begin{lemma}
For $\z \in \mathcal D \mathfrak M^\infty(U)$ the following hold:

{\rm (a)} \cite[Proposition 3.1]{CheFri} $|\divi \z| \ll \mathcal H^{N-1}$.
	
{\rm (b)} \cite[Proposition 3.4]{ACM}, \cite[Lemma 2.4]{GMP} $\divi \z$ can be decomposed as
	\[\divi \z = \divi^a \z + \divi^j \z + \divi^c \z, \ \text{with} \ \left\{ \begin{array}{l}  \divi^a \z \ll \mathcal L^N, \smallskip \\   \divi^c \z(B) = 0, \text{if } \mathcal H^{N-1}(B) < \infty, \end{array}\right.\]
and
\begin{equation} \label{jump-div}
\divi^j \z = \big([\z, \nu_u]^+ -[\z, \nu_u]^-\big)\mathcal H^{N-1}\res J_u.
\end{equation}
\end{lemma}

\subsection{Pairing measures}

 Let $U \subset \R^N$ be  any  open set.

\begin{definition}[Anzellotti's pairing] Given the pair $\z \in \mathcal D \mathfrak M^\infty_{\rm loc}(U)$ and $u \in BV_{\rm loc}(U) \cap L^\infty_{\rm loc}(U)$, we define  $(z, Du) \in \mathcal D'(U)$ by means of
\[\big<(\z, Du), \varphi\big>:= -\int_U u^\ast \, \varphi \, d (\divi \z)  - \int_U u \, \z \cdot \nabla \varphi \, dx \quad \text{for all} \ \varphi \in C_0^\infty(U).\]
\end{definition}
\noindent This is a well-defined Radon measure on $U$, and $(\z, Du) \ll |Du|$ (see \cite[Lemma 5.1]{Cas}). Moreover, we can write (see \cite[Theorem 3.2]{CheFri})
\[(\z, Du) = (\z, Du)^a + (\z, Du)^j+ (\z, Du)^c, \quad \text{with } (\z, Du)^a = \z \cdot \nabla u  \, \mathcal L^N,\]
and a further characterization for the jump part (cf. \cite[Theorem 4.12]{CraCi}) is
\begin{equation} \label{pairing-jump}
(\z, Du)^j = \frac{[\z, \nu_u]^+ + [\z, \nu_u]^-}{2} \, (u^+-u^-) \mathcal H^{N-1}\res J_u.
\end{equation}

\begin{lemma}[Generalized Green's formula and product rules] \label{prod-rule-div} For a vector field $\z \in \mathcal D \mathfrak M^\infty_{\rm loc}(U)$ and $u \in BV_{\rm loc}(U) \cap L^\infty_{\rm loc}(U)$, it holds
\begin{equation} \label{green}
\int_{\Omega} u^\ast d(\divi \, \z)+ \int_\Omega d (\z,Du) = \int_{\partial \Omega} [\z, \nu] u \, d \mathcal H^{N-1},
\end{equation}
 for any open and bounded $\Omega \subset U$ with Lipschitz continuous boundary.
Moreover, we get the product rules:

{\rm (a)} ${\rm div}(u \z) = u^\ast {\rm div} \, \z + (\z,Du)$  as measures,

{\rm (b)} \!\cite[Proposition 4.9]{CraCi} For $\z \in \mathcal D \mathfrak M^\infty_{\rm loc}(U)$, and $u, v \in BV_{\rm loc}(U)\cap L^\infty_{\rm loc}(U)$, 
\begin{equation} \label{prule-pair}
(v\z,Du) = v^\ast (\z,Du) + \frac{(u^+ - u^-)(v^+ - v^-)}{4} \divi \z \res (J_u \cap J_v).
\end{equation}
\end{lemma}

 Hereafter, unless otherwise specified, we will always identify any $BV$ function $v$ with its precise representative $v^\ast$, without explicit use of  $\ast$.

\subsection{Anisotropic total variation}

For an anisotropy $F$,  an open set $U$ and $u \in BV(U)$, the $F$-total variation or anisotropic total variation $F(Du)= |Du|_{_F}\in\mathfrak{M}(U)$ of $Du$ is given by
\begin{align} \label{anis-decomp}
|Du|_{_F} &:=  F(\nu_u)|Du|\\
&= F(\nabla u) \mathcal L^N + F(\nu_u) |u^+-u^-|  \, \mathcal H^{N-1}\res J_u + F(\nu_u) |D^c u|.
\end{align}
We will use the notation $|Dv|_{_F}(U) = \displaystyle \int_U dF(Dv)$.

Here are the main properties:
\begin{lemma} \label{prop-ATV} Let $U\subset\R^N$ be a bounded open set with Lipschitz continuous boundary and $u\in BV(U)$. Then,

{\rm (a)}  the $F$-total variation of $Du$ in an open bounded set $A\subset U$ with Lipschitz continuous boundary can be computed as follows (cf. \cite{AmBe})
\begin{align} \label{FtotV}
\hspace*{-0.9cm}|Du|_{_F}(A)
 =\sup\!\bigg\{\!\int_A \!\! u \, \divi \z \, :    \z \in C^1_c(A;\R^N),  \, F^\circ(\z) \leq 1 \text{ a.e.\ in } A\bigg\}.
\end{align}

{\rm (b)} $\displaystyle \int_U d F(Du)$ is $L^1(U)$-lower semicontinuous on $BV(U)$.
\end{lemma}

 Notice that (a) is proven in \cite[p.185]{moll_05} while (b) is obtained in \cite[Theorem 5.1]{AmBe}. On the other hand, we note that if $\divi\z\in L^2(U)$, the following result has been shown in \cite[Proposition 4.1]{CFM} for item (a) and (b) can be found in \cite[proof of Theorem 7]{moll_05}. However, we need to work out the technical details to be able to obtain the same conclusions under milder assumptions.

\begin{lemma}  \label{prop-ATV-2}
Let $U$ be an open set with bounded Lipschitz boundary, and $\z\in \mathcal D\mathfrak{M}^\infty_{\rm loc}(U)$. Then the following Cauchy-Schwarz type inequalities hold:

{\rm (a)} $\big|(\z, Du)\big| \leq \|F^\circ(\z)\|_{L^\infty(U)} F(Du)$, for all  $u\in BV_{\rm loc}(U)\cap L^\infty_{\rm loc}(U)$.

{\rm (b)} $|[\z,\nu]|\leq \|F^\circ(\z)\|_{L^\infty(U)} F(\nu)$, $\mathcal H^{N-1}$-a.e. on $\partial U$. 
\end{lemma}

\begin{proof}
In order to prove (a), we first extend $u$ to all of $\R^N$; let $g\in W^{1,1}(\R^N\setminus U)$ such that $g=u$ at $\partial U$ (see, for example, \cite[Theorem 18.18]{Leoni}).
  Then, let $\overline u:=u\chi_U+g\chi_{\R^N\setminus U}$ and consider $u_n:=\overline u\star \rho_n \in C^\infty(\R^N)$, with $\{\rho_n\}_{n\in \N}$ a sequence of mollifiers. We have that $u_n\to u$ in $L^1(U)$ and, by Lemma \ref{conv-u-ast}, $u_n\to u^*$ $\mathcal H^{N-1}$-a.e. as $n\to \infty$. Moreover, by \cite[Proposition 3.7]{Ambrosio},
  $$\int_A |\nabla u_n|\xrightarrow[]{n\to \infty}|Du|(A) \ \  \text{ for every open set } A\Subset U {\rm \ with \ }|Du|(\partial A)=0.$$
  Then, by Resethnyak's continuity theorem (\cite[Theorem 2.39]{Ambrosio}) and \eqref{anis-decomp},
  $$\int_A F(\nabla u_n)\xrightarrow[]{n\to \infty}|Du|_{_F}(A) \ \  \text{ for every open set } A\Subset U {\rm \ with \ }|Du|(\partial A)=0.$$
    Therefore, given an open set $A\Subset U$ such that $|Du|(\partial A)=0$ and $\varphi\in C_0^\infty(A)$,  by the definition of Anzellotti's pairing,
  \begin{align}\langle (\z, Du),\varphi\rangle= &-\int_A   u \z\cdot\nabla \varphi-\int_A u^*\varphi \, d(\divi\z)\\
  	= &- \lim_{n\to \infty} \bigg(\int_A   u_n \z\cdot\nabla \varphi+\int_A u_n\varphi \, d(\divi\z)\bigg)\\
  =&\lim_{n\to \infty} \int_A \varphi \, \z\cdot\nabla u_n\leq \|F^\circ (\z)\|_{L^\infty(U)}\|\varphi\|_{L^\infty(A)}\lim_{n\to \infty}\int_A F(\nabla u_n)\\
  =& \|F^\circ (\z)\|_{L^\infty(U)}\|\varphi\|_{L^\infty(A)} |Du|_{_F}(A),\end{align}
  which implies (a).

  Before proving (b), we need the following approximation result:

{\noindent\bf Claim:} Given $\z\in \mathcal D\mathfrak{M}^\infty_{\rm loc}(U)$ there exists $\z_n\in C^\infty(\overline U;\R^N)$ such that $$\z_n\to \z {\rm \ in \ } L^1_{\rm loc}(U;\R^N),\quad \divi \z_n\stackrel{*}\rightharpoonup \divi \z {\rm \ in \ } \mathfrak{M}(U).$$
  The claim is proved as the corresponding result for vector fields with divergence in $L^p(U)$; see \cite[Lemma 2.3]{kohn-temam} (cf.\cite[Lemma 3.6]{Andreu-Caselles-Mazon-ibero}).

  Take a partition of unity of $\overline U$, $\{\theta_j\}_{j=1}^k\subset C^\infty_0(\R^N)$, $0\leq\theta_j\leq 1$, satisfying the following: if ${\rm supp \ }\theta_{j}\cap \partial U\neq\emptyset$, then there is a bounded open cone $C_{j}$ with vertex $0$ such that $x+C_{j}\cap\overline U=\emptyset$ for all $x\in\partial U \cap {\rm supp \ }\theta_{j}$ and there exists $r>0$ such that $x-C_{j}\subseteq U$ for all $x\in\partial U\cap({\rm supp \ }\theta_{j}+B_r(0))$. Moreover, for each $j\in \{1,\ldots, k\}$, take a mollifier $\rho_j\in C^\infty(\R^N)$ with the property that, if ${\rm supp \ }\theta_{j}\cap \partial U\neq\emptyset$, then ${\rm supp \ }\rho_{j}\subset C_{j}$.

 Letting $\rho_{j,n}(x):=n^N\rho_j(nx)$, we define
   $$\z_{n}:=\sum_{j=1}^k \rho_{j,n}\star(\theta_j\overline \z)\in C^\infty(\R^N;\R^N), \quad \text{where} \quad \overline \z = \left\{ \begin{array}{ll} z, & \text{in } U \smallskip \\ 0, & \text{in } \R^N \setminus U.\end{array}\right.$$

   Note that, if $x\in\partial U\cap {\rm supp \ }\theta_{j}$, then $x-C_j\subset U$ and ${\rm supp \ }\rho_j\subset C_j$ thus 
   \begin{align}
   (\rho_{j,n}\star(\theta_j\overline \z))(x)&=\int_{x-C_j}\rho_{j,n}(x-y)\theta_j(y)\overline \z(y)dy=\int_{U}\rho_{j,n}(x-y)\theta_j(y) \z(y)dy.
   \end{align}

 \vspace*{-0.2cm}

By construction, $\z_n\to\z$ in $L^1_{\rm loc}(U;\R^N)$. On the other hand, for any $j=1,\ldots, k$,
  \begin{equation} \label{div-decomp}
  \divi (\theta_j\overline\z)=\overline{\divi(\theta_j\z)} +[\theta_j\overline\z,\nu]\mathcal H^{N-1}\res_{\partial U},
  \end{equation}
  where $\overline{\divi(\theta_j\z)}(E)=\divi(\theta_j\z)(E\cap U)$ for any $E\subset \R^N$. Now, for $x\in\R^N$,
  \begin{align}
  \rho_{j,n}\star\left([\theta_j\overline\z,\nu]\mathcal H^{N-1}\res_{\partial U}\right)(x)&=\int_{\partial U\cap \, {\rm supp \ }\theta_j} \rho_{j,n}(x-y)[\theta_j\overline\z,\nu](y)d\mathcal H^{N-1}(y),
  \end{align}
  thus $\rho_{j,n}\star\left([\theta_j\overline\z,\nu]\mathcal H^{N-1}\res_{\partial U}\right)=0$ if $\partial U\cap \, {\rm supp \ }\theta_j=\emptyset$. Otherwise, one gets $y+C_{j}\subseteq \R^N\setminus\overline U$ for all $y\in\partial U\cap {\rm supp \ }\theta_{j}$ and  ${\rm supp \ }\rho_{j}\subset C_{j}$. In this case, if $y\in\partial U\cap \, {\rm supp \ }\theta_j$, $\rho_{j,n}(x-y)\neq 0$ only if $x\in y+{\rm supp \ }\rho_j\subset \R^N\setminus\overline U$ so
  \begin{align}
  \int_{\partial U\cap \, {\rm supp \ }\theta_j} \rho_{j,n}(x-y)[\theta_j\overline\z,\nu](y)d\mathcal H^{N-1}(y)=0 \ \ \text{ for every } x\in \overline U,
  \end{align}
  that is, $\rho_{j,n}\star\left([\theta_j\overline\z,\nu]\mathcal H^{N-1}\res_{\partial U}\right)=0$ in $\overline U$. Hence, using \eqref{div-decomp}, we can write
$$(\divi(\rho_{j,n}\star(\theta_j\overline\z)))\res_U=(\rho_{j,n}\star \divi (\theta_j\overline \z))\res_U=\big(\rho_{j,n}\star \overline{\divi (\theta_j\z)}\big)\res_U.$$
Consequently,
\vspace*{-0.2cm}
    \[\divi(\z_n)\res_U=\sum_{j=1}^k \divi\left(\rho_{j,n}\star(\theta_j\overline \z)\right)\res_U=\sum_{j=1}^k \big(\rho_{j,n}\star\overline{\divi (\theta_j \z)}\big)\res_U.\]
    Letting $n\to\infty$, since $\big(\overline{\divi (\theta_j \z)}\big)\res_U=\divi (\theta_j \z)$, we finally obtain the claim.

 Therefore, by the definition of the normal trace $[\z,\nu]$ (cf \cite[Section 3]{ACM})  and the classical Green's formula, we obtain, for any $0\leq \Psi\in C^\infty(U)\cap W^{1,1}(U)$ with bounded support,
\begin{align}
\int_{\partial U} \Psi\ \z_n\cdot\nu &=\int_{\partial U}\Psi \ \sum_{j=1}^k (\rho_{j,n}\star(\theta_j\overline \z))\cdot\nu \\
&=\int_{U}\Psi \ \sum_{j=1}^k \rho_{j,n}\star\overline{\divi (\theta_j \z)}+\int_{U}\nabla\Psi \cdot \sum_{j=1}^k (\rho_{j,n}\star(\theta_j \z))\\
& \xrightarrow[]{n\to \infty} \int_{U}\Psi \ d(\divi \z)+\int_{U}\nabla\Psi \cdot \z=\int_{\partial U} \Psi \ [\z,\nu].
\end{align}
  Then, by approximation we have that
  $$\lim_{n\to\infty}\int_{\partial U} \Psi\ \z_n\cdot\nu= \int_{\partial U} \Psi \ [\z,\nu]$$
  for every $\Psi\in W^{1,1}(U)$. Consequently, for any such $\Psi$ and by \eqref{CS-ineq},
  \begin{align*}\left|\int_{\partial U} \Psi \ [\z,\nu]\right|
   \leq &\lim_{n\to\infty} \sum_{j=1}^k\!\int_{\partial U}\Psi(x)\left(\int_U \rho_{j,n}(x-y)\theta_j(y)\left|\z(y)\cdot\nu(x)\right|{\rm d}y\right) \!{\rm d }x\!\\
  \leq & \lim_{n\to\infty}\sum_{j=1}^k\!\int_{\partial U}\!\!\Psi(x)F(\nu(x))\!\left(\int_U \!\!\rho_{j,n}(x-y)\theta_j(y)F^\circ(\z(y)){\rm d}y\right) \!{\rm d }x\!\\ \leq & \|F^\circ(\z)\|_{L^\infty(U)}\int_{\partial U}\Psi F(\nu).\end{align*}
  This last inequality directly gives (b).
\end{proof}

 \section{Definition and uniqueness of solutions for A-IMCF} \label{sec-def-flow}

\subsection{Definition of solutions and first properties} \label{def-sol1}
As in \cite{MazSer13, MazSer15} we introduce the following notion of solution for \eqref{EquationVp1}:
\begin{definition} \label{weak-IMCF-def}
A  non-negative $v \in BV_{\rm loc}(\Omega^e)\cap L^\infty_{\rm loc}(\Omega^e)$ is said to be a weak solution of \eqref{EquationVp1} if there is a vector field ${\rm {\bf z}} \in \mathcal D \mathfrak M^\infty_{\rm loc}(\Omega^e)$, with $F^\circ({\rm {\bf z}}) \leq 1$, satisfying
\begin{equation} \label{weak-IMCF}
{\rm div} \,\z = (\z, Dv) = F(Dv) \quad  \text{as measures in } \ \Omega^e,
\end{equation}
together with the boundary conditions
\begin{equation} \label{IMCF-bdry}
\left\{\begin{array}{ll}
v = 0 & {\mathcal H}^{N-1}\text{-a.e. on } \, \partial \Omega,
\medskip \\  v(x) \to \infty & \text{as} \ |x| \to \infty.
\end{array}\right.
\end{equation}
\end{definition}

\smallskip

 \begin{remark}
In the terminology by Moser (see e.g. \cite[Definition 1.1]{Moser15}), the behaviour at infinity in \eqref{IMCF-bdry} means that all our solutions are {\it proper}, which geometrically implies that solutions stay bounded at finite times, and hence the corresponding sublevel sets are precompact. Thus we do not need to add the latter as an extra assumption as in \cite{HuIl}. As pointed out in \cite{KoNi}, there are some manifolds that do not admit a proper solution, even in more regular settings. We stress that the lack of properness leads to non-uniqueness of solutions (cf.~\cite[Proposition 2.4]{Moser08}).
\end{remark}

\begin{lemma} \label{no-jump}
A weak solution $v$ of \eqref{EquationVp1} in the sense of Definition \ref{weak-IMCF-def} has no jump part, i.e., it satisfies $D^j v = 0$.
\end{lemma}

\begin{proof} Let $\mu\in \mathfrak{M}(\Omega^e)$ be such that $\mu\ll\mathcal H^{N-1}\res J_v$ and let us denote by $\theta^{J_v}[\mu]$ the Radon-Nikodým derivative of $\mu$ with respect to $\mathcal H^{N-1}\res J_v$. Now by \eqref{weak-IMCF} and \eqref{pairing-jump} we can write
\begin{align*}
\theta^{J_v}[F(Dv)] & =\theta^{J_v}[(\z, Dv)] = \frac{[\z, \nu_v]^+ + [\z, \nu_v]^-}{2} \, (v^+-v^-) \\ & = \Big(\frac{\theta^{J_v}[\divi \z]}{2} + [\z, \nu_v]^- \Big) (v^+-v^-), 
\end{align*}
which follows from \eqref{jump-div}. Using now that, by \eqref{weak-IMCF}, $\divi \z = F(Dv)$ as measures, we get  by  means of \ref{prop-ATV-2} (b) 
\begin{align*}
\theta^{J_v}[F(Dv)] \bigg(\!1-\frac{v^+-v^-}{2}\!\bigg) &= [\z, \nu_v]^- (v^+-v^-) \leq \|F^\circ(\z)\|_{L^\infty} F(\nu_v) |v^+-v^-|
\\ & \leq F(\nu_v) |v^+-v^-|  = \theta^{J_v}[F(Dv)],
\end{align*}	
where we applied $F^\circ(\z) \leq 1$ and \eqref{anis-decomp}. This leads to $v^+ \geq v^-$ $\mathcal{H}^{N-1}$-a.e. Redoing the above argument, but substituting $ [\z, \nu_v]^-$ instead of  $ [\z, \nu_v]^+$ by means of \eqref{pairing-jump}, one can reach the reverse inequality. Thus $v^+ = v^-$ $\mathcal{H}^{N-1}$-a.e., and by \eqref{Du-dec} the claim follows.
\end{proof}

In \cite{GMP}, $BV$ functions with null singular set with respect to the Hausdorff measure $\mathcal H^{N-1}$ were named as Diffuse $BV$ functions and denoted by $DBV$. We recall some of their properties, which  will be useful in the sequel.
\begin{lemma}\label{DBV}
  Let $U\subset\R^N$ be an open set, $v\in DBV_{\rm loc}(U)$ and $f:\R\to\R$ be a locally Lipschitz function. Then, the  following chain rules hold:
  \begin{itemize}
    \item[{\rm (a)}] $D(f \circ v) = f'(v) Dv.$ \smallskip
\item[{\rm (b)}] For any $\z\in \mathcal D \mathfrak M^\infty_{\rm loc}(U)$,   if $f$ is non-decreasing, we have\smallskip
\begin{itemize}
\item[{\rm (b1)}] $(\z,D (f\circ v))=f'(v)(\z,D v).$\smallskip
\item[{\rm (b2)}] If additionally $(\z,Dv)=F(Dv)$,  then \[(\z,D(f\circ v))=F(D(f\circ v)).\]
  \end{itemize}
\end{itemize}
\end{lemma}
\begin{proof}
  (a) follows from the standard chain rule formula for BV functions \cite[Theorem 3.99]{Ambrosio}.
  In order to show (b1), we proceed as follows:
 \begin{align*}(\z,D(f\circ v))=&\frac{(\z,D(f\circ v))}{|D(f\circ v)|}|D(f\circ v)|=\frac{(\z,Dv)}{|Dv|}|D(f\circ v)|\\ =& f'(v)\frac{(\z,Dv)}{|Dv|}|Dv|=f'(v)(\z,Dv),\end{align*}
  where the second equality is due to \cite[Proposition 4.5 (iii)]{CraCi} and the third one comes from (a). Now, we obtain (b2) by combining the previous items in the statement with \eqref{anis-decomp}. Indeed,
 \begin{align*}(\z,D(f\circ v))= & f'(v)(\z,Dv)=f'(v)F(Dv)\\ =& f'(v)F(\nu_v)|Dv|=F(\nu_v)|D(f\circ v)|=|F(D(f\circ v))|.\end{align*}
 Here we have used that $\nu_{f\circ v}=\frac{D (f\circ v)}{|D (f\circ v)|}=\nu_v$, $D(f\circ v)$-a.e.
\end{proof}

\begin{corollary} \label{div=0}
	Let $v$ be a weak solution of \eqref{weak-IMCF} with associated vector field $\z$. Then
	the measure $\divi \z = F(Dv)$ vanishes when restricted to the level sets, that is,
	\[\divi \z\res \{v = t\}= 0 \quad \text{for all} \quad t \geq 0.\]
\end{corollary}

\begin{proof}
As in \cite[Remark 4.4]{MazSer15} one proves that $Dv\res \{v = t\}= 0$ for any $v \in DBV_{\rm loc}(\Omega^e)\cap L^\infty_{\rm loc}(\Omega^e)$, and then use the anisotropic total variation and equation \eqref{weak-IMCF}. 	
\end{proof}

\subsection{Uniqueness of the anisotropic IMCF}

\begin{theorem}\label{unique}
  Suppose that $u_1$ and $u_2$ are a non-negative weak supersolution and a non-negative weak subsolution of \eqref{EquationVp1} in the sense that $u_1$, $u_2 \in  DBV_{\rm loc}(\Omega^e)\cap L^\infty_{\rm loc}(\Omega^e)$ and there exist vector fields $\z_i \in \mathcal D \mathfrak M^\infty_{\rm loc}(\Omega^e)$, with $F^\circ(\z_i) \leq 1$,  $i=1,2$, satisfying
$${\rm div} (\z_1) \leq (\z_1, Du_1) = F(Du_1) \quad  \text{as measures in } \ \Omega^e,$$
and
$${\rm div} (\z_2) \geq (\z_2, Du_2)= F(Du_2) \quad  \text{as measures in } \ \Omega^e,$$
together with boundary conditions
$$
\left\{\begin{array}{ll}
u_i = \varphi_i & {\mathcal H}^{N-1}\text{-a.e. on } \, \partial \Omega,
\medskip \\  u_i(x) \to \infty & \text{as} \ |x| \to \infty
\end{array}\right.,
$$  with $\varphi_i\in L^1(\partial\Omega)$, such that
$\varphi_2\le \varphi_1$ $\mathcal{H}^{N-1}$-a.e. Then $u_2\le u_1$ a.e. in $\Omega^e$.
\end{theorem}
\begin{proof}
  Let $T^h:\R\to (-\infty,h]$ be the truncation operator defined by
\begin{equation} \label{trunc-up}
T^h(s):=\left\{\begin{array}{ll}
                   s & \hbox{if } s\le h, \\
                   h & \hbox{if } s> h.
                 \end{array}\right.
\end{equation}

  Suppose by contradiction that $u_2> u_1$ holds on a set of positive measure. Then, {since $u_1\ge 0$}, there exists $h>0$ such that $T^h(u_2)> u_1$ on a set of positive measure. Therefore, we can find $\varepsilon>0$ such that $l:=\esssup \left\{T^h(u_2)-\frac{u_1}{1-\varepsilon}\right\}-\varepsilon>0$ (note that $l\le h$). Define $\overline u_1:=\frac{u_1}{1-\varepsilon}$ and $\overline u_2:=T^h(u_2)-l$.  We have
  \begin{equation}\label{supersolution}
  {\rm div}({\rm {\bf z}}_1)\leq(1-\varepsilon) F(D\overline u_1),\end{equation}
    \begin{equation}\label{subsolution}
{\rm div}({\rm {\bf z}}_2)\geq F(D u_2),\end{equation}
  \begin{equation}\label{esssup}
  \esssup \{\overline u_2-\overline u_1\}=\varepsilon,\end{equation}
  and we note that $\overline u_2(x)\geq \overline u_1(x)$ implies $u_2(x)>u_1(x)$.

 Hereafter, we use the notation $u_+:=\max\{u,0\}$ for the positive part of a function $u$. Since $\overline u_1\to\infty$ as $|x|\to\infty$ and $\overline{u}_2\le h$, we have $(\overline u_2-\overline u_1)_+=0$ on $\partial B_R$ for $R$ large enough. Moreover, $(\overline u_2-\overline u_1)_+\le (u_2- u_1)_+=0$ on $\partial \Omega$. Therefore, if we integrate $(\overline u_2-\overline u_1)_+$ over $B_R\setminus\overline\Omega$ with respect to the measures in \eqref{subsolution} and \eqref{supersolution}, take the difference between them and integrate by parts we get
  \begin{align}\label{eq:unique}&-\int_{B_R \setminus \overline\Omega} d({\rm {\bf z}}_2-{\rm {\bf z}}_1,D(\overline u_2-\overline u_1)_+)\\
  &\geq\int_{B_R \setminus \overline\Omega} (\overline u_2-\overline u_1)_+ \, d (F(D u_2)-F(D\overline u_1))+\varepsilon\int_{B_R \setminus \overline\Omega} (\overline u_2-\overline u_1)_+ \, d F(D\overline u_1).\end{align}
  We observe that, by Lemma \ref{DBV} (b2), $({\rm {\bf z}}_i,D\overline u_i)=F(D\overline u_i)$, $i=1,2$ . Therefore, since $F^\circ({\rm {\bf z}}_i)\leq 1$, $i=1,2$, by \eqref{CS-ineq} and Lemma \ref{DBV} (b1)   the first integral on the left hand side is nonpositive:
  $$\int_{B_R \setminus \overline\Omega} d({\rm {\bf z}}_2-{\rm {\bf z}}_1,D(\overline u_2-\overline u_1)_+)=\int_{B_R \setminus \overline\Omega}\chi_{\{\overline u_2-\overline u_1\geq 0\}} \,d({\rm {\bf z}}_2-{\rm {\bf z}}_1,D(\overline u_2-\overline u_1))\geq 0.$$
 On the other hand, multiplying the inequality of measures ${\rm div} ({\rm {\bf z}}_1)\leq F(D\overline u_1)$ by $\frac12((\overline u_2-\overline u_1)_+)^2$, integrating by parts and using again Lemma \ref{DBV} (b1), we obtain
    \begin{align}\label{claim}
    \int_{B_R \setminus \overline\Omega} (\overline u_2-\overline u_1)_+ \, d (F(D\overline u_2)-F(D\overline u_1))\geq  -\int_{B_R \setminus \overline\Omega} \frac{((\overline u_2-\overline u_1)_+)^2}{2} \, d F(D\overline u_1).
  \end{align}
  Therefore, since
  $$\int_{B_R \setminus \overline\Omega} (\overline u_2-\overline u_1)_+ \, d (F(D u_2)-F(D\overline u_1))\geq \int_{B_R \setminus \overline\Omega} (\overline u_2-\overline u_1)_+ \, d (F(D\overline u_2)-F(D\overline u_1)),$$
  \eqref{claim} together with \eqref{eq:unique} yield
  \begin{align}\int_{B_R \setminus \overline\Omega} (\overline u_2-\overline u_1)_+\left(\varepsilon-\frac{(\overline u_2-\overline u_1)_+}{2}\right) d F(D\overline u_1)\leq 0.\end{align}
  Now, letting $R\to\infty$ we have
    $$\int_{\Omega^e} (\overline u_2-\overline u_1)_+\left(\varepsilon-\frac{(\overline u_2-\overline u_1)_+}{2}\right) d F(D\overline u_1)\leq 0$$
  thus, by \eqref{esssup},
  $$\int_{\Omega^e} (\overline u_2-\overline u_1)_+  \, d F(D\overline u_1)=0.$$
  This, in turn, yields
  $$\int_{\Omega^e} (\overline u_2-\overline u_1)_+  \, d |D\overline u_1|=0.$$
  Moreover, from \eqref{eq:unique}, we get
  $$\int_{\Omega^e} (\overline u_2-\overline u_1)_+  \, d |D \overline u_2|\le \int_{\Omega^e} (\overline u_2-\overline u_1)_+  \, d |D u_2|=0.$$
  Finally,
  \begin{align}
    \frac12|D (( \overline u_2-\overline u_1)_+)^2|(\Omega^e)=\int_{\Omega^e} (\overline u_2-\overline u_1)_+  \, d |D(\overline u_2-\overline u_1)| \\
    \leq  \int_{\Omega^e} (\overline u_2-\overline u_1)_+d |D\overline u_1|+ \int_{\Omega^e} (\overline u_2-\overline u_1)_+  \, d |D\overline u_2|=0,
  \end{align}
  which implies that $(\overline u_2-\overline u_1)_+$ is equivalent to a constant in $\Omega^e$. Since $\overline u_2-\overline u_1\to 0$ as $|x|\to\infty$, such a constant is zero. This is in contradiction with \eqref{esssup}.
\end{proof}

\section{Definition and uniqueness of solutions of the approximating problems} \label{Moser-aprox}

\subsection{Anisotropic $p$-Laplace equation and change of variable} \label{def-sol-2}
 \begin{definition}
  A function $u\in L^{pN/(N-p)}(\Omega^e)$ such that $\nabla u \in L^p(\Omega^e)$ is said to be a weak solution of the problem \eqref{EquationU} if its trace on $\partial \Omega$ is equal to 1,  $\displaystyle \lim_{|x| \to \infty} u(x) = 0$ and there exists $\z \in \partial F(\nabla u)$ such that
  \[{\rm div}\big( F^{p-1}(\nabla u) \z\big) = 0\]
 in the weak sense, meaning that
  \begin{equation} \label{weakU}
  \int_{\Omega^e} F^{p-1}(\nabla u) \z \cdot \nabla \varphi = 0 \quad \text{for every} \quad \varphi \in W^{1,p}_0(\Omega^e).
  \end{equation}
 \end{definition}

In \cite{CRMS} we proved an existence theorem for weak solutions of the above problem:
\begin{theorem} \label{CRMS-thm}
If $1 < p < N$ and $\Omega \subset \R^N$ is a bounded domain with Lipschitz-continuous boundary, there exists a weak solution $u$ of the problem \eqref{EquationU} satisfying the estimates
\begin{equation}\label{boundsu-U}\left(F^\circ\right)^{\frac{p-N}{p-1}}\bigg(\frac{x}{r_1}\bigg)\le u(x) \le \left(F^\circ\right)^{\frac{p-N}{p-1}}\bigg(\frac{x}{r_2}\bigg) \quad \hbox{a.e. in } \Omega^e
\end{equation}
for any constants  $0<r_1<r_2$ such that $\mathcal W_{r_1}\subseteq \Omega\subseteq \mathcal W_{r_2}$.
\end{theorem}

We will see that the transformation \eqref{change-var} gives a weak solution to the problem \eqref{EquationV} in the following sense:
\begin{definition} \label{def:sol:V}
We say that a function $v \in W^{1,p}_{\rm loc}(\Omega^e)$ is a {\it weak solution} of \eqref{EquationV} if its trace on $ \partial\Omega$ is 0, $\displaystyle \lim_{|x| \to \infty} v(x) = \infty$ and there exists $\tilde \z \in \partial F(\nabla v)$ such that
\[{\rm div}\big( F^{p-1}(\nabla v) \tilde \z\big) = F^p(\nabla v) \quad \text{in the weak sense} ,\]
that is,
\begin{equation} \label{test:vp}
-\int_{\Omega^e} F^{p-1}(\nabla v) \tilde \z \cdot \nabla \varphi = \int_{\Omega^e} F^{p}(\nabla v) \varphi,
\end{equation}
for every  $\varphi\in L^\infty(\Omega^e)$ with compact support whose distributional gradient  $\nabla \varphi$ belongs to $L^p(\Omega^e; \R^N)$.
\end{definition}

\begin{proposition}\label{variable-change}
	 If $u$ is a weak solution of \eqref{EquationU}, then $v$ given by \eqref{change-var} is a weak solution of \eqref{EquationV}.
\end{proposition}

\begin{proof}

Notice that
\begin{equation} \label{zVsz}
\z \in \partial F(\nabla u) \quad \text{ if and only if } \quad \tilde \z = -\z \in  \partial F(\nabla v).
\end{equation}
 Indeed, by 1-homogeneity of $F$, we have
\[\tilde \z \cdot \nabla v = (p-1) \z\cdot \frac{\nabla u}{u} =\frac{p-1}{u} F(\nabla u) = F\Big(\frac{1-p}{u} \nabla u\Big) = F(\nabla v).\]

Because of the lower bound in \eqref{boundsu-U}, for each open and bounded set $B \subset \Omega^e$, we can find a constant $C = C(p, r_1, N, B) > 0$ so that $u \geq C$ on $B$. As $r \mapsto 1/r^{p-1}$ is Lipschitz-continuous away from zero, a well-known result by Stampacchia \cite{Stam} ensures that $1/u^{p-1} \in W^{1,p}(B)$, and hence for any $\varphi \in  W^{1,p}_0(B)\cap L^\infty(B)$ it holds that
\[\Big(\frac{p-1}{u}\Big)^{p-1} \varphi \in W^{1,p}_0(B) \cap L^\infty(B)\]
can be used as a test function in \eqref{weakU}. Accordingly, for $\z \in \partial F(\nabla u)$ we have
\[0 =  \int_B  F^{p-1}\Big(\frac{p-1}{u} \nabla u\Big)  \z \cdot \nabla \varphi - \int_B   F^p\Big(\frac{p-1}{u} \nabla u\Big) \varphi,\]
where we applied the homogeneity of $F$ and that $\nabla u \cdot \z = F(\nabla u)$. Now, as $\tilde \z = -\z \in \partial F(\nabla v)$, this can be rewritten as
\[0 =  -\int_B  F^{p-1}(\nabla v)  \tilde \z \cdot \nabla \varphi - \int_B   F^p( \nabla v) \varphi,\]
as desired.
\end{proof}

As an immediate consequence of \autoref{CRMS-thm}, we conclude:
\begin{corollary} \label{exist-vp}
If  $\Omega \subset \R^N$ is a bounded domain with Lipschitz-continuous boundary, there exists a weak solution $v$ of the problem \eqref{EquationV}. Moreover, $v$ satisfies the estimates
\begin{equation}\label{boundsu} (N-p) \log\bigg(F^\circ\left(\frac{x}{r_2}\right)\bigg)\le v(x) \le (N-p) \log\bigg(F^\circ\left(\frac{x}{r_1}\right)\bigg) \quad \hbox{a.e. in } \Omega^e.
\end{equation}
for any constants  $0<r_1<r_2$ such that $\mathcal W_{r_1}\subseteq \Omega\subseteq \mathcal W_{r_2}$.	
\end{corollary}

\subsection{Comparison theorems for sub- and supersolutions}

 Consider the following generalization of  \eqref{EquationV} with generic Dirichlet boundary condition:
  \begin{equation}\label{EquationVp1varphi}
 \left\{ \begin{array}{ll} {\rm div}\big(F^{p-1}(\nabla v) \partial F(\nabla v)\big) = F^p(\nabla v) & \text{in }  \Omega^e \smallskip \\
 v = \varphi & \text{on }  \partial\Omega \smallskip \\
 v \to  \infty & \text{as }  |x| \to \infty \end{array} \right..
 \end{equation}

 \begin{theorem}\label{comparison_vp}
 Suppose that $v_1$ and $v_2$ are a non-negative weak supersolution and a non-negative weak subsolution of \eqref{EquationVp1varphi} with boundary values $\varphi_1\ge \varphi_2\in W^{1-\frac1{p}, p}(\partial\Omega)$, respectively, on $\partial\Omega$, i.e., $v_i \in W^{1,p}_{\rm loc}(\Omega^e)$, its trace on $\partial \Omega$ is equal to $\varphi_i$,  $\displaystyle \lim_{|x| \to \infty} v_i(x) = \infty$ and there exists $\z_i \in \partial F(\nabla u_i)$ such that
\[{\rm div}\big( F^{p-1}(\nabla v_1) \z_1\big) \leq F^p(\nabla v_1) \quad \text{in the weak sense}, \]
and
\[{\rm div}\big( F^{p-1}(\nabla v_2) \z_2\big) \geq F^p(\nabla v_2) \quad \text{in the weak sense} .\]
Then $v_1\ge v_2$ a.e. in $\Omega^e$.
 \end{theorem}
 \begin{proof}
  The proof follows the ideas of the proof of Theorem \ref{unique}, suitably adapted to this more regular case.
 Let $T^h:\R\to (-\infty,h]$ be the truncation operator defined by \eqref{trunc-up}.

 We argue by contradiction. Suppose that $v_2>v_1$ holds on a set of positive measure. Then, since $v_1\ge 0$, there exists $h>0$ such that $T^h(v_2)> v_1$ on a set of positive measure. Therefore, we can find $\varepsilon>0$ such that $l:=\esssup \left\{T^h(v_2)-\frac{v_1}{1-\varepsilon}\right\}-\frac{\varepsilon}{p}>0$ (note that $l\le h$). Define $\overline v_1:=\frac{v_1}{1-\varepsilon}$ and $\overline v_2:=T^h(v_2)-l$. We have
  \begin{equation}\label{supersolutionPp}
  {\rm div}\big( F^{p-1}(\nabla \overline v_1)\z_1\big) \leq (1-\varepsilon)F^p(\nabla \overline v_1) \quad \text{in the weak sense},\end{equation}
    \begin{equation}\label{subsolutionPp}
{\rm div}\big( F^{p-1}(\nabla \overline v_2)\z_2\big) \geq F^p(\nabla \overline v_2) \quad \text{in the weak sense},\end{equation}
  \begin{equation}\label{esssupPp}
  \esssup \{\overline v_2-\overline v_1\}=\frac{\varepsilon}{p}.\end{equation}

  Since $v_1\to \infty$ as $|x|\to\infty$, we have $\overline v_1>h$ on $\partial B_R$ for all $R$ large enough. Therefore, since $\overline v_2 < h$, $(\overline v_2-\overline v_1)_+=0$ on $\partial B_R$  for all $R$ large enough. Moreover, $(\overline v_2-\overline v_1)_+\le (v_2- v_1)_+=0$ on $\partial \Omega$. Accordingly, if we multiply \eqref{subsolutionPp} and \eqref{supersolutionPp} by $(\overline v_2-\overline v_1)_+$ over $B_R\setminus\overline\Omega$, take the difference between them and integrate by parts, we get
  \begin{align}\label{eq:uniquePp}&-\int_{B_R \setminus \overline\Omega} \left(F^{p-1}(\nabla \overline v_2)\z_2-F^{p-1}(\nabla \overline v_1)\z_1\right)\cdot\nabla(\overline v_2-\overline v_1)_+\\
  &\geq\int_{B_R \setminus \overline\Omega} (F^p(\nabla \overline v_2)-F^p(\nabla\overline v_1))(\overline v_2-\overline v_1)_++\varepsilon\int_{B_R \setminus \overline\Omega} F^p(\nabla\overline u_1)(\overline v_2-\overline v_1)_+.\end{align}
  We recall that by \eqref{charFp} it holds $pF^{p-1}(\nabla \overline v_i) \z_i\in \partial F^p(\nabla\overline v_i)$, $i=1,2$, by the chain rule for subdifferentials \eqref{charFp}. Therefore, by the monotonicity of $\partial F^p$, the integral on the left hand side is nonpositive. On the other hand, multiplying the inequality
  $${\rm div}\big( F^{p-1}(\nabla \overline v_1)\z_1\big)\leq F^p(\nabla\overline v_1)$$
  by $\frac{p}{2}((\overline v_2-\overline v_1)_+)^2$ and integrating by parts, we obtain
    \begin{align}\label{claimPp}
    \int_{B_R \setminus \overline\Omega} (F^p(\nabla\overline v_2)-F^p(\nabla\overline v_1))(\overline v_2-\overline v_1)_+ \geq  -\int_{B_R \setminus \overline\Omega}F^p(\nabla\overline v_1)\frac{p((\overline v_2-\overline v_1)_+)^2}{2}.
  \end{align}
  Therefore, \eqref{claimPp} together with \eqref{eq:uniquePp} yield
  \begin{align}\int_{B_R \setminus \overline\Omega} (\overline v_2-\overline v_1)_+\left(\varepsilon-\frac{p(\overline v_2-\overline v_1)_+}{2}\right) F^p(\nabla \overline v_1)\leq 0.\end{align}
  Now, letting $R\to\infty$ we have
    $$\int_{\Omega^e} (\overline v_2-\overline v_1)_+\left(\varepsilon-\frac{p(\overline v_2-\overline v_1)_+}{2}\right) F^p(\nabla \overline v_1)\leq 0.$$
  thus, by \eqref{esssupPp},
  $$(\overline v_2-\overline v_1)_+  F(\nabla \overline v_1)=0 \quad \text{a.e. in } \Omega^e.$$
  Moreover, from \eqref{eq:uniquePp}, we get
  $$(\overline v_2-\overline v_1)_+  F(\nabla \overline v_2)=0 \quad \text{a.e. in } \Omega^e.$$
  Finally,
  \begin{align*}
    \frac12 \left|\nabla (( \overline v_2-\overline v_1)_+)^2\right|& =\left|(\overline v_2-\overline v_1)_+\nabla(\overline v_2-\overline v_1)\right|\leq (\overline v_2-\overline v_1)_+ \left(|\nabla \overline v_1|+|\nabla \overline v_2| \right)\\
    &\leq \frac1c(\overline v_2-\overline v_1)_+ \left(F(\nabla \overline v_1)+F(\nabla \overline v_2) \right)=0\quad \text{a.e. in } \Omega^e.
  \end{align*}
  This implies that $(\overline v_2-\overline v_1)_+$ is equivalent to a constant in $\Omega^e$. Since $(\overline v_2-\overline v_1)_+=0$ for $|x|$ large enough, this constant has to be zero, which contradicts with \eqref{esssupPp}.
 \end{proof}

\begin{remark}
  We actually do not need to assume $\lim_{|x| \to \infty} v_2(x) = \infty$.
 \end{remark}

 \subsection{Uniqueness of the $p$-capacitary potentials}	
 As a by-product of the techniques in the previous proof, we manage to achieve uniqueness of solutions of \eqref{EquationU}, which enhances our existence and regularity results in \cite{CRMS}. More precisely, we show that the transformation \eqref{change-var} allows us to obtain a comparison principle 
  which yields uniqueness of solutions of this problem. In fact, setting $p^* = \frac{p N}{N-p}$,

\begin{theorem}\label{UniquenessPcapacity}
  Suppose that $u_1$ and $u_2$ are a weak subsolution and a weak supersolution of \eqref{EquationUvarphi} with $D = \Omega^e$ that have boundary values $\varphi_1\le \varphi_2\in W^{1-\frac1{p}, p}(\partial\Omega)$, respectively, on $\partial\Omega$, i.e., $u_i\in L^{p^*}(\Omega^e)$,  $\nabla u_i \in L^p(\Omega^e)$, its trace on $\partial \Omega$ is equal to $\varphi_i$,  $\displaystyle \lim_{|x| \to \infty} u_i(x) = 0$ and there exists $\z_i \in \partial F(\nabla u_i)$ such that
  \[{\rm div}\big( F^{p-1}(\nabla u_2) \z_2\big) \leq 0 \leq {\rm div}\big( F^{p-1}(\nabla u_1) \z_1\big)\]
 in the weak sense. Then $u_1\le u_2$ a.e. in $\Omega^e$.
\end{theorem}
\begin{proof}
Let $T_{h_1}^{h_2}:\R\to [h,+\infty)$ be the truncation operator defined by
$$T_{h_1}^{h_2}(s):=\left\{\begin{array}{ll}
                   h_1 & \hbox{if } s\le h_1, \smallskip \\
                   s & \hbox{if } h_1\le s\le h_2, \smallskip\\
                   h_2 & \hbox{if } s\ge h_2.
                 \end{array}\right.$$

  Suppose by contradiction that $u_1> u_2$ holds on a set of positive measure. Then, there exist $0<h<1$ and $h<h'$ such that either $T_{h}^{h'}(u_1)> T_{h}^{h'}(u_2)$ or $T_{h}^{h'}(-u_1)<T_{h}^{h'}(-u_2)$ holds on a set of positive measure. We may assume without loss of generality that the former holds, since otherwise we can consider the solutions $\tilde u_i=-u_i$, $i=1,2$, of \eqref{EquationU} and get to a contradiction in the same way. Moreover, we can assume that $h'=1$, since otherwise we consider the solutions $\tilde u_i=\frac{1}{h'}u_i$, $i=1,2$.

  Now, let $v_i=v_i^h:=(1-p)\log (T_h^1(u_i))\in W^{1,p}_{\rm loc}(\Omega^e)$, $i=1,2$. Note that $0\le v_i\le (1-p)\log h$ and $v_1\ge v_2$ on $\partial \Omega$. We claim that,  similarly as in Proposition \ref{variable-change}, one can show the inequalities
  \[{\rm div}\big( F^{p-1}(\nabla v_1)(-\z_1)\big) \leq F^p(\nabla v_1) \quad \text{in the weak sense},\]
  and
    \[{\rm div}\big( F^{p-1}(\nabla v_2)(-\z_2)\big) \geq F^p(\nabla v_2) \quad \text{in the weak sense}.\]

   From this point on, the proof follows exactly the one of Theorem \ref{comparison_vp}. Therefore, we are led to a contradiction (we omit the details).
\end{proof}

 Recall that the solutions of \eqref{EquationUvarphi} with $D = \Omega^e$ and $\varphi \equiv 1$ are known to be minimizers of the anisotropic $p$-capacity (see \cite{CRMS} for more details). Accordingly,

\begin{corollary}
There exists a unique minimizer  in $\{u\in L^{p^\ast}(\R^N)\,:\, \nabla u \in L^p(\R^N), \, u\geq 1 \text{ in } \Omega\}$ of ${\rm Cap}^F_p(\overline\Omega)$  defined as in \eqref{pcap-def}, that is, $p$-capacitary functions are unique.
\end{corollary}

With a similar proof to that of Theorem \ref{UniquenessPcapacity} we obtain a comparison principle, thus uniqueness, of weak solutions of \eqref{EquationUvarphi} on bounded domains.
\begin{theorem}
  Suppose that $u_1$ and $u_2$ are two weak solutions of \eqref{EquationUvarphi} with $D = \Omega$ and boundary values $\varphi_1\le \varphi_2\in W^{1-\frac1{p}, p}(\partial\Omega)$, respectively. Then $u_1\le u_2$ a.e. in $\Omega$.
\end{theorem}

 Summing up, the results in this subsection give altogether a proof of Corollary \ref{uniq-CRMS1}.

\section{Existence of solutions to A-IMCF} \label{existence-sect}

\begin{theorem} \label{exist-thm}
		If $\Omega \subset \R^N$ is a bounded domain with Lipschitz-continuous boundary, there exists a weak solution $v$ of the problem \eqref{EquationVp1}. Moreover, $v$ satisfies the estimates
		\begin{equation} \label{barrier-lim}
		(N-1) \log\bigg(F^\circ\left(\frac{x}{r_2}\right)\bigg)\le v(x) \le (N-1) \log\bigg(F^\circ\left(\frac{x}{r_1}\right)\bigg) \quad \hbox{a.e. in } \Omega^e.
		\end{equation}
		for any constants  $0<r_1<r_2$ such that $\mathcal W_{r_1}\subseteq \Omega\subseteq \mathcal W_{r_2}$.	
\end{theorem}

\begin{proof}
By Corollary \ref{exist-vp}, we know that for each $p>1$ there exists a weak solution $v_p$ of \eqref{EquationV} in the sense of Definition \ref{def:sol:V}. Moreover, by uniqueness, $v_p = (1-p) \log u_p$, where $u_p$ is the unique weak solution of \eqref{EquationU}.

In particular, \eqref{test:vp} holds for some $\z_p \in \partial F(\nabla v_p)$. We aim to construct a solution of problem \eqref{weak-IMCF} as a subsequential limit of $\{v_p\}_{p >1}$ as $p\searrow 1$. We divide the proof into several steps.

{\bf Step 1: BV-estimate.} Let $B \subset \Omega^e$ be an open and bounded set and $0 \leq \varphi \in C_c^\infty(\Omega^e)$ such that $\varphi \equiv 1$ in $B$. Then we take $\varphi^p$ as a test function in \eqref{test:vp} and, by means of  Young's inequality  with $p':= \frac{p}{p-1}$, we get
\begin{align*}
\int_{\Omega^e}   F^p( \nabla v_p) \varphi^p  &=  -p \int_{\Omega^e}  F^{p-1}(\nabla v_p) \varphi^{p-1}   \z_p \cdot \nabla \varphi \\ &  \leq p \int_{\Omega^e}  \big(F(\nabla v_p) \varphi\big)^{p-1}   F^\circ(\z_p) F(\nabla \varphi)
\\ & \leq \frac{p}{p'} \int_{\Omega^e}  F^{p}(\nabla v_p) \varphi^{p} + \int_{\Omega^e} F^p(\nabla \varphi).
\end{align*}
Here we have used \eqref{CS-ineq} and $F^\circ(\z_p) \leq 1$, where the latter follows from \eqref{char-subdif-1hom}. Consequently, for all $1 < p < 3/2$, we have
\begin{equation} \label{gra-Lp-vp}
\int_{B}   |\nabla v_p|^p \leq \left(\frac{C}{c}\right)^p \frac1{2-p} \int_{{\rm supp}\, \varphi} |\nabla \varphi|^p \leq 2 \left(\frac{C}{c}\right)^{\frac3{2}} \int_{{\rm supp}\, \varphi} (1 +|\nabla \varphi|^{\frac{3}{2}}) \leq \mathcal C,
\end{equation}
where $c, C$ are the constants  from \eqref{FequivEuclid}. Hereafter $\mathcal C$ denotes any constant which depends on $B$, $c$ and $C$, but is independent of $p$, whose concrete meaning may change from line to line. Now H\"older's inequality leads to
\[\int_{B}   |\nabla v_p| \leq \mathcal C \qquad \text{for all} \quad 1<p<3/2.\]

Since, by \eqref{boundsu}, $\{v_p\}$ is pointwise bounded independently of $p$ for $p \in(1,3/2)$, we conclude that $\{v_p\}$ is bounded in $W^{1,1}(B)$ for every $B \subset \Omega^e$ open and bounded, and hence by the embedding theorem \cite[Corollary 3.49]{Ambrosio}, we get a function $v \in BV_{\rm loc}(\Omega^e)$ such that, up to subsequences and as $p\searrow 1$,
\begin{align}
&v_p  \to v \quad \text{in} \ L^q_{\rm loc}(\Omega^e) \quad \text{for every} \ 1 \leq q <\frac{N}{N-1} \quad\text{and a.e. in } \ \Omega^e, \label{conv-vp-v} \smallskip \\
&\nabla v_p \res B \rightharpoonup Dv\res B \ \text{weakly-$\ast$ as measures for each open bounded } B \subset \Omega^e.
\end{align}

Notice that the behaviour of $v$ as $|x|\to \infty$, as well as the estimates in \eqref{barrier-lim}, follow directly by taking limits in \eqref{boundsu}.

{\bf Step 2: Existence of $\z$.} The aim is to find $\z \in \mathcal D \mathfrak M^\infty_{\rm loc}(\Omega^e)$ with $F^\circ(\z) \leq 1$  and satisfying
\begin{equation} \label{conv-z}
	F^{p-1}(\nabla v_p)  \z_p \rightharpoonup \z \quad \text{weakly in } L^1_{\rm loc}(\Omega^e; \R^N) \quad \text{as} \ p\searrow 1.
\end{equation}
 First, for $B \Subset \Omega^e$ fixed, we show that the sequence $\{F^{p-1}(\nabla v_p)  \z_p\}_{1 < p <3/2}$ is equibounded in $L^1(B; \R^N)$. Indeed, by \eqref{FequivEuclid}, \eqref{F0equivEuclid}, Hölder's inequality and \eqref{gra-Lp-vp}, we have for $1 \leq q < \frac{p}{p-1}$
\begin{align*}
\int_B \big|F^{p-1}(\nabla v_p) \z_p\big|^q &  \leq C^q\int_B F^{q(p-1)}(\nabla v_p)F^\circ(\z_p)^q \leq C^{p q} \int_B |\nabla v_p|^{q(p-1)}  \\ & \leq C^{pq} \|\nabla v_p\|_{L^p(B)}^{q(p-1)} \,|B|^{\alpha(p)} \leq \mathcal C
\end{align*}
with $\alpha(p) \to 1$ as $p \to 1$. Therefore, for any $1 \leq q< \infty$, we conclude that there exists a subsequence so that
\begin{equation}
	F^{p-1}(\nabla v_p)  \z_p \rightharpoonup \z_{B} \quad \text{weakly in } L^q(B; \R^N).
\end{equation}

{\bf Claim 1.} {\it $\z_{_B} \in L^\infty(B; \R^N)$ with $F^\circ(\z_{B})\leq 1$. }
\begin{proof}[Proof of Claim 1]
We adapt to our setting an argument from the proof of \cite[Lemma 1]{ABCM}. For any $k > 0$, define
\[B_{p,k} =\big\{x \in B \ : \ F\big(\nabla v_p(x)\big) > k\big\}.\]
By means of \eqref{gra-Lp-vp} and \eqref{FequivEuclid}, we can estimate the measure of this set by
\begin{equation} \label{meas-B-bdd}
\big|B_{p,k}\big| \leq \int_{B_{p,k}} \frac{F^p\big(\nabla v_p(x)\big)}{k^p} \leq \frac{\mathcal C}{k^p} \quad \text{for every }   1 < p < 3/2 \text{ and } k>0.
\end{equation}
Now,
\[F^{p-1}(\nabla v_p) \z_p \chi_{B_{p,k}} \rightharpoonup \g_k:=\z_{B}\,\chi_{B_{p,k}} \quad \text{ weakly in } L^1(B; \R^n) \text{ as } p\searrow 1.\]
Choose $\Phi \in L^\infty(B; \R^N)$ with $F(\Phi) \leq 1$. Then, since $F^\circ(\z_p) \leq 1$ because $\z_p \in \partial F(\nabla v_p)$, Cauchy-Schwarz \eqref{CS-ineq} and H\"older's inequality lead to
\[\bigg|\int_{B_{p,k}} F^{p-1}(\nabla v_p) \z_p \cdot \Phi\bigg| \leq  C^{p-1} \big\|\nabla v_p\big\|_{L^p(B)}^{p-1} \,\big|B_{p,k}\big|^{1/p} \leq \frac{\mathcal C}{k},\]
where we used \eqref{gra-Lp-vp} and \eqref{meas-B-bdd}. Taking limits as $p \searrow 1$ and choosing $\Phi \propto {\rm sgn}(\g_k)$, we get
\[\int_{B} |\g_k | \leq  \frac{\mathcal C}{k}, \quad \text{for every} \ k > 0.\]

On the other hand, we have that
\[F^\circ\big(F^{p-1}(\nabla v_p) \z_p \, \chi_{B \setminus B_{p,k}}\big) =  F^{p-1}(\nabla v_p)\chi_{B \setminus B_{p,k}} F^\circ(\z_p) \leq   k^{p-1}\] for any $p >1$. Accordingly,
\[F^{p-1}(\nabla v_p) \z_p \, \chi_{B \setminus B_{p,k}} \rightharpoonup \bm{f}_k :=\z_{B}\,\chi_{B \setminus B_{p,k}} \quad \text{ weakly in } L^1(B; \R^N) \text{ as } p\searrow 1,\]
with $F^\circ({\bm f}_k) \leq 1$. In short, for any $k >0$, we can write ${\z_{_B}} = \g_k + \bm{f}_k$ so that
\[F^\circ({\z_{_B}})\leq F^\circ(\g_k) + F^\circ(\bm{f}_k) \leq 1,\]
from which we conclude that $F^\circ({\z_{_B}}) \leq 1$, as desired.
\end{proof}

Next, we choose an increasing sequence $\{B_n\}_{n\in\mathbb N}$ with $B_n \Subset \Omega^e$ such that $\displaystyle\cup_{n\in \mathbb N} B_n = \Omega^e$. One can perform a diagonal argument to establish that (up to a subsequence) \eqref{conv-z} holds,
where $\z: \Omega^e \rightarrow \R^N$ is a vector field such that $\z|_B = \z_B$ for each $B\Subset \Omega^e$ and, consequently, $\z \in L^\infty(\Omega^e;\R^N)$ with $F^\circ(\z) \leq 1$ by virtue of Claim 1.

It remains to check that $\z \in \mathcal D \mathfrak M^\infty_{\rm loc}(\Omega^e)$. Now, \eqref{test:vp} and \eqref{gra-Lp-vp} ensure that  the sequence $\big\{\divi(F^{p-1}(\nabla v_p) \z_p)\big\}$ is bounded in $L^1_{\rm loc}(\Omega^e; \R^N)$. Therefore, we can find a Radon measure $\mu$ such that, as $p\searrow 1$, we have the subsequential convergence
\[{\rm div}(F^{p-1}(\nabla v_p)  \z_p) \rightharpoonup \mu \quad \text{weakly-}\ast  \text{as measures}.\]
In particular, recalling also \eqref{conv-z}, we have that, for any open and bounded set $B \subset \Omega^e$ and $\varphi\in C^1_c(B)$,
\[\int_B F^{p-1}(\nabla v_p)  \z_p \cdot \nabla \varphi \to  - \int_{B} \varphi \, d \mu = \int_{B} {\z} \cdot \nabla \varphi,\]
thus $\mu =\divi \z$. The latter implies that $\divi \z$ has locally bounded total variation, as desired.

{\bf Step 3.}  {\bf $\bm{{\rm div}\, \z = (\z,Dv)=F(Dv)}$ as measures.}
Take $0\leq \varphi \in C^\infty_c(\Omega^e)$ as a test function in \eqref{test:vp} to have
\begin{equation} \label{step2-p}
\int_{\Omega^e} F^{p-1}(\nabla v_p)  \z_p \cdot \nabla \varphi + \int_{\Omega^e} F^{p}(\nabla v_p) \varphi = 0.
\end{equation}

On the other hand, Young's inequality yields
\begin{align}
\int_{\Omega^e}  F(\nabla v_p) \, \varphi &\leq \frac1{p} \int_{\Omega^e} F^p(\nabla v_p) \varphi + \frac1{p'} \int_{\Omega^e} \varphi \label{aux-Young}
\\ & = - \frac1{p}\int_{\Omega^e} F^{p-1}(\nabla v_p)  \z_p \cdot \nabla \varphi + \frac{p-1}{p} \int_{\Omega^e} \varphi. \nonumber
\end{align}
This inequality joint with lower semicontinuity (see Lemma \ref{prop-ATV} (b)) give
\begin{align*}
\int_{\Omega^e} \varphi \, d F(Dv) &\leq \liminf_{p\searrow 1} \int_{\Omega^e} F(\nabla v_p)\varphi \leq - \int_{\Omega^e} \z \cdot \nabla \varphi,
\end{align*}
which follows from \eqref{conv-z}. Hence,
\begin{equation} \label{easy-ineq3}
-{\rm div} \, \z + F(Dv) \leq 0.
\end{equation}

To prove the reverse inequality, let us take $\varphi \in C^\infty_c(\Omega^e)$ and use $e^{-v_p} \varphi$ as a test function in \eqref{test:vp}. By means of \eqref{char-subdif-1hom}, we reach
\[0 = \int_{\Omega^e} e^{-v_p} F^{p-1}(\nabla v_p) \z_p \cdot \nabla \varphi.\]
Since $e^{-v_p} \leq 1$, by \eqref{conv-vp-v} and the dominated convergence theorem we have that $e^{-v_p}$ converges to  $e^{-v}$ in $L^q({\rm supp} \,\varphi)$ for all  $1 \leq q < \infty$. Then we may take limits as $p\searrow 1$, and get
\[0 = \int_{\Omega^e} e^{-v} \z \cdot \nabla \varphi,\]
where we applied \eqref{conv-z}. In particular, we have
\begin{equation} \label{div-exp0}
\divi(e^{-v} \z) = 0 \quad \text{in } \mathcal D'(\Omega^e).
\end{equation}

Here we can use the product rule from Lemma  \ref{prod-rule-div} (a) and write
\begin{align}\label{ineqexpz}
0 =  e^{-v}\, \divi \z + \big(\z,D(e^{-v})\big) \geq e^{-v} \, F(Dv) - F(D e^{-v}),
\end{align}
where the inequality follows from \eqref{easy-ineq3}, the Cauchy-Schwarz type inequality from Lemma \ref{prop-ATV-2} (a) and $F^\circ(\z) \leq 1$. By definition of the precise representative (cf.\ Lemma \ref{conv-u-ast}) and of $F(Dv)$ in \eqref{anis-decomp}, we can write
\begin{equation} \label{base_ineq}
0 \geq  e^{-v} F(Dv)^d + \frac{e^{-v^+} + e^{-v^-}}{2} F(\nu_v) |v^+-v^-|\mathcal H^{N-1}\res J_v- F(D e^{-v}).
\end{equation}

Notice that, as we are integrating over a compact subset by the choice of $\varphi$, we can apply the chain rule for BV functions (see \cite[Theorem 3.99]{Ambrosio}) and the 1-homogeneity of $F$ to reach
\begin{align*}
F(D e^{-v}) = e^{- v} F(Dv)^d + \big|e^{-v^+} - e^{-v^-}\big| F(\nu_v) \mathcal H^{N-1}\res J_v.
\end{align*}
Substitution of the latter into \eqref{base_ineq} yields
\begin{align}\label{convexjump}
 0 \geq  \bigg(\frac{e^{-v^+} + e^{-v^-}}{2} |v^+-v^-|- \big|e^{-v^+} - e^{-v^-}\big| \bigg)F(\nu_v) \mathcal H^{N-1}\res J_v.
\end{align}
Given $a\geq 0$, we let 
$$f(h) =e^{-a}\left( \frac{h}{2}(e^{-h} + 1)+ e^{-h}-1\right).$$
It holds that $f(h)\geq 0$ in $\R^+$ and $f(h)=0$ if, and only if, $h=0$. Therefore, taking $a=\min\{v^-(x),v^+(x)\}$ and $h=|v^+-v^-|$, having in mind \eqref{convexjump}, we obtain that $v^-=v^+$ $\mathcal H^{N-1}-$a.e in $J_v$; i.e. $v\in DBV_{\rm loc}(\Omega^e)$. Consequently,  from \eqref{div-exp0} and Lemma \ref{DBV} (b1), we get
$$0 = \divi(-e^{-v} \z) = -e^{-v}\divi \z +e^{-v}(\z,Dv) \quad {\rm in \ }\mathcal D'(\Omega^e),$$ which implies, by means of \eqref{easy-ineq3} and Lemma \ref{prop-ATV-2} (a), that  $$F(Dv)\leq \divi \z=(\z,Dv)\leq F(Dv) \quad {\rm in \ }\mathcal D'(\Omega^e).$$

Finally, an approximation argument ensures 
$$\divi\z= {(\z, Dv) = F(Dv)},\quad {\rm as \ measures.}$$

{\bf Step 4. Boundary condition.} The aim is to show that $v = 0$ on $\partial \Omega$ in the sense of traces. Take $R >0$ big enough so that $\Omega \subset B_R(0)$ and take $\zeta v_p$ as a test function in \eqref{test:vp}, where $\zeta \in C_c^\infty(\R^N)$ given by
\[0\leq \zeta \leq 1 \qquad \zeta(x) =\left\{ \begin{array}{ll} 1, & \text{if } |x|\leq R \smallskip \\ 0,  & \text{if } |x|\geq 2R. \end{array} \right.\]

 Doing so, we obtain
 $$\int_{\Omega^e}  F^p(\nabla v_p)v_p\zeta = -\int_{\Omega^e} F^p(\nabla v_p)\zeta - \int_{\Omega^e}   F^{p-1}(\nabla v_p)  v_p \, \z_p\cdot \nabla \zeta. $$
From here, by means of Young's inequality, we reach
\begin{align} \label{step4-aux}
\int_{\Omega^e} \zeta F(\nabla v_p) &+ \int_{\Omega^e} \! F(\nabla v_p) v_p \zeta \nonumber \\ &\leq \frac{1}{p}\int_{\Omega^e}  \zeta F^p(\nabla v_p) + \frac{1}{p'} \int_{\Omega^e}  \zeta + \frac{1}{p} \int_{\Omega^e}  v_p \zeta F^p(\nabla v_p) + \frac{1}{p'} \int_{\Omega^e}  v_p \zeta \nonumber
\\ &= \frac{p-1}{p} \int_{\Omega^e}  (1+ v_p) \zeta -\frac1{p} \int_{\Omega^e}  v_p F^{p-1}(\nabla v_p)  \z_p \cdot \nabla \zeta.
\end{align}

Set $\Omega_R:= B_{2R} \setminus \overline{\Omega}$ and recall that $v_p$ vanishes on $\partial \Omega$. Then, by lower semicontinuity (see \cite[Theorem 5 and 7]{moll_05}), we can write
 $$\int_{\Omega_R} \zeta dF(Dv)+ \int_{\partial \Omega} v F(\nu_v) \le\liminf_{p \searrow 1}\int_{\Omega_R} \zeta F(\nabla v_p),$$
as well as
$$\int_{\Omega_R} \frac{\zeta}{2} \, dF\left(Dv^2\right)+ \int_{\partial \Omega} \frac{v^2}{2}F(\nu_v) \le\liminf_{p \searrow 1}\int_{\Omega_R} \zeta v_p F(\nabla v_p).$$
Adding these inequalities up and taking limits in \eqref{step4-aux}, we have
\begin{align} \label{step4-aux2}
\int_{\Omega_R} \!\!\zeta dF(Dv) + \int_{\partial \Omega} \!\! F(\nu_v) v + \int_{\Omega_R} \!\! \zeta v \, dF(Dv)  + \int_{\partial \Omega} \!\! F(\nu_v)  \frac{v^2}{2} \leq -\int_{\Omega_R} \!\! v \, \z \cdot \nabla \zeta,
\end{align}
where we also applied \eqref{conv-z}.

On the other hand, if we combine the product rule from Lemma \ref{prod-rule-div} and the generalized Green's formula \eqref{green}, we reach
\begin{align*}
-\int_{\Omega_R} \!\! v \, \z \cdot \nabla \zeta & = \int_{\Omega_R} \zeta v   \, d(\divi \z) + \int_{\Omega_R} \zeta \, d (\z, Dv) -\int_{\partial \Omega} v [\z, \nu_v] \\ & = \int_{\Omega_R} \zeta v   \, dF(Dv) + \int_{\Omega_R} \zeta d F(Dv) -\int_{\partial \Omega} v [\z, \nu_v],
\end{align*}
which follows from the equalities in Step 3.

Accordingly, by substitution of the latter into \eqref{step4-aux2}, all the integrals on $\Omega_R$ cancel out, and we conclude
\[\int_{\partial \Omega} \!\! \big(F(\nu_v)  v +  v [\z, \nu_v]\big) \,   + \frac1{2} \int_{\partial \Omega} \!\!  F(\nu_v)  v^2 \leq 0.\]
Notice that $v [\z, \nu_v] \geq -v F(\nu_v) \|F^\circ(\z)\|_\infty \geq -v F(\nu_v)$ by Lemma \ref{prop-ATV-2}(b), and hence on the left hand side we have the sum of two non-negative quantities. Thus we get that $v$ vanishes on $\partial \Omega$ in the sense of traces, as desired.
\end{proof}

\section{Extra regularity for domains satisfying a uniform interior ball condition} \label{cont-harnack}

 We recall that, if $\Omega$ satisfies the $\mathcal W_r$-condition stated in Definition \ref{UIBC}, then solutions to the approximating problems \eqref{EquationU} are Lipschitz continuous (\cite[Theorem 1.5]{CRMS}). Accordingly, it is natural to wonder if for the limiting problem \eqref{EquationVp1} we also get additional regularity in case that one can slide a Wulff shape of uniform size along the boundary of $\Omega$.

 In this section we show that, under these circumstances, the solution to \eqref{EquationVp1} is actually continuous and satisfies a Harnack inequality. The strategy to prove this closely follows that of the corresponding statements for the Euclidean IMCF with obstacles (\cite[Theorem 2.5 and Proposition 2.7]{Moser08}).

  	\begin{theorem} \label{cont-thm}
 Let $r>0$ and suppose that $\Omega$ satisfies the $\mathcal W_r$-condition. Then the unique solution of \eqref{EquationVp1} is continuous.
 	\end{theorem}
 \begin{proof}
 Let $x_0\in\partial\Omega$. Then, there exists $y_0\in \R^N$ such that
 \[\mathcal W_r+y_0\subseteq \overline{\Omega} \quad \text{ and } \quad x_0\in\partial \left(\mathcal W_r+y_0\right).\]
 Define
 $$\eta_p(x):=\left(\frac{F^\circ(x-y_0)}{r}\right)^{\frac{p-N}{p-1}}.$$
 Then $\eta_p$ is a solution of \eqref{EquationUvarphi} with $\eta_p\leq 1$ on $\partial\Omega$. Therefore, by Theorem \ref{UniquenessPcapacity}, $\eta_p\le u_p$ in $\Omega^e$. In particular,
 $$v_p(x)=(1-p)\log u_p(x) \le (N-p)\log \left(\frac{F^\circ(x-y_0)}{r}\right) \quad \text{in } \Omega^e.$$
 Then, for any $\delta>0$ and $x\in  (\mathcal W_\delta+x_0)\cap\Omega^e$, we have
 \begin{align}
v_p(x)&\le (N-p)\log \left(\frac{F^\circ(x-x_0)}{r}+\frac{F^\circ(x_0-y_0)}{r}\right)\\
&\leq (N-p)\log \left(\frac{\delta}{r}+1\right).
 \end{align}
 Consequently, for any $\varepsilon>0$ we may choose $\delta>0$ small enough, depending  on $r, \varepsilon$ and $N$ (but independent of $x_0$), such that
 $$v_p(x)\le \varepsilon \quad \text{for } x\in (\mathcal W_\delta+x_0)\cap\Omega^e.$$

For $|{\rm \bf  h}|<\delta$, define
\[v_p^{{\rm \bf  h}}(x):=v_p(x + {\rm \bf  h})+\varepsilon \quad \text{ for } \quad x\in\Omega-{\rm \bf  h}.\]
 Then (see Figure \ref{fig:moser}),
 \[v_p^{\rm \bf  h}=\varepsilon \ \text{ on } \ \partial (\Omega-{\rm \bf  h}) \qquad \text{ and } \qquad v_p^{\rm \bf  h}\geq \varepsilon \ \text{ on } \ \partial \Omega \cap (\Omega^e-{\rm \bf  h}).\]
 \vspace*{-0.5cm}
 \begin{figure}[H]
 	\includegraphics[scale=0.5]{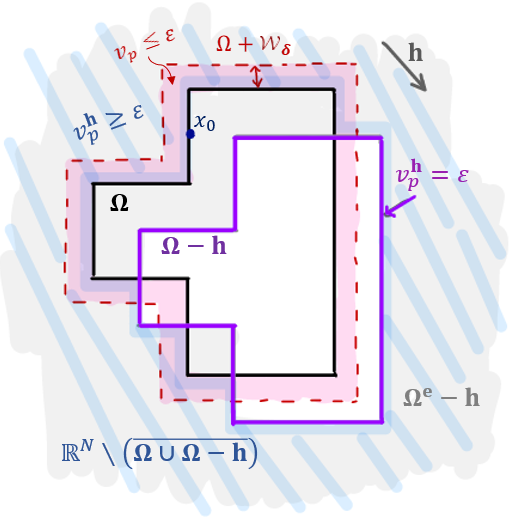}\vspace*{-0.3cm}\caption{Different regions for the comparison arguments}\label{fig:moser}
 	\end{figure}
 \vspace*{-0.5cm}
 Now, $v_p$ and $v_p^{\rm \bf  h}$ are two solutions of \eqref{EquationVp1varphi} on the exterior domain $\R^N\setminus \overline{E}$, with $E=\Omega\cup (\Omega-{\rm \bf  h})$,  so that $v_p\le v_p^{\rm \bf  h}$ on $\partial E$. Therefore, by Theorem \ref{comparison_vp} we get that
 \[v_p\le v_p^{\rm \bf  h} \quad \text{ on } \quad \R^N\setminus \overline E\supset \R^N\setminus \overline{(\Omega+ \tfrac1{c}\mathcal{W}_\delta)}.\]

Since the same is true for $-{\rm \bf  h}$ instead of ${\rm \bf  h}$, we conclude that $|v_p(x+{\rm \bf  h})-v_p(x)|\le \varepsilon$ for every $|{\rm \bf  h}|<\delta$ and $x\in \R^N\setminus \overline{(\Omega+ \tfrac1{c}\mathcal{W}_\delta)}$. That is, the $v_p$ are equicontinuous outside of any neighbourhood of $\Omega$. In particular, they are equicontinuous on any compact subset of $\Omega^e$ thus the Arzelà-Ascoli Theorem ensures that $v$ is continuous on $\Omega^e$.
 \end{proof}

\begin{theorem}\label{Harnack}
Let $v$ the unique solution to \eqref{EquationVp1} for a given domain $\Omega$ which satisfies a $\mathcal W_s$ condition for some $s >0$. There exists $\mathcal C>0$ such that, if $x_0\in\Omega^e$ and $\rho>0$ satisfy $\mathcal W_{2 \rho}(x_0)\subset \Omega^e$, then $${\rm osc}_{\mathcal W_\rho(x_0)} v\leq \mathcal C.$$
\end{theorem}

\begin{proof}
Assume that $p\in \left(1,\frac32\right]$. Let  $u_p\in W^{1,\infty}(\Omega^e)$ 
be a weak solution of \eqref{EquationU}. Given $x_0\in \Omega^e$, let $\rho>0$ be such that $\mathcal W_{2\rho}(x_0)\subset \Omega^e$. Since the equation in \eqref{EquationU} is invariant under translations and dilations, we may suppose that $x_0=0$ and $\rho=1$.

For $q\in \R\setminus \{0\}$ and $r\in(0,2]$ we denote
$$I(q,r):=\left(\fint_{\mathcal W_r} u_p^q\right)^{\frac1q}.$$
Now, by Hölder's inequality, if $q_1<q_2$,
$$I(q_1,r)\leq \frac{1}{|\mathcal W_r|^\frac{1}{q_1}}\left(\int_{\mathcal W_r}u_p^{q_2}\right)^\frac{1}{q_2}|\mathcal W_r|^\frac{q_2-q_1}{q_1q_2}=I(q_2,r),$$
thus $I(q,r)$ is non-decreasing in $q$. Moreover, since $u_p$ is continuous,
$$I(\infty,r)=\lim_{q\to\infty}I(q,r)=\sup_{\mathcal W_r}u_p$$
and
$$I(-\infty,r)=\lim_{q\to-\infty}I(q,r)=\inf_{\mathcal W_r}u_p.$$

Let $v_p:=(1-p)\log u_p$. Following the proof of Proposition \ref{variable-change}, we get that $v_p$ satisfies
$${\rm div}\big(F^{p-1}(\nabla v_p) \partial F(\nabla v_p)\big)= F^p(\nabla v_p) \quad \text{in } \mathcal W_2$$
in the weak sense; that is, there exists $\z\in\partial F(\nabla v_p)$ such that
$$-\int_{\mathcal W_2} F^{p-1}(\nabla v_p) \z \cdot \nabla \varphi = \int_{\mathcal W_2} F^{p}(\nabla v_p) \varphi,
$$
for every  $\varphi\in L^\infty(\mathcal W_2)$ with compact support  in $\mathcal W_2$ whose distributional gradient  $\nabla \varphi$ belongs to $L^p(\mathcal W_2; \R^N)$.
Now, for any such $\varphi$, by \eqref{CS-ineq} and Hölder's inequality we have
\begin{align}
\int_{\mathcal W_2} F^{p}(\nabla v_p) \varphi^p &= -\int_{\mathcal W_2} p\varphi^{p-1}F^{p-1}(\nabla v_p)\z\cdot\nabla \varphi\\
&\leq \left(\int_{\mathcal W_2} F^{p}(\nabla v_p) \varphi^p\right)^{\frac{p-1}{p}}\left(\int_{\mathcal W_2}p^pF^p(\nabla\varphi)\right)^\frac{1}{p}.
\end{align}
Thus
$$\int_{\mathcal W_2} F^{p}(\nabla v_p) \varphi^p \leq p^p\int_{\mathcal W_2}F^p(\nabla\varphi).$$
In particular, for any $y_0\in {\mathcal W_2}$ and $r$ such that $\mathcal W_{2r}(y_0)\subseteq \mathcal W_2$, we may appropriately choose $\varphi$ to obtain
$$\fint_{\mathcal W_r(y_0)} F^{p}(\nabla v_p) \leq \frac{2^N p^p}{r^p}.$$
Then, the Poincaré inequality for  anisotropic balls (proved as in \cite[Section 5.8 Theorem 2]{Evans}) yields
$$\fint_{\mathcal W_r(y_0)}\left|v_p(x)- \overline{v_p}^r\right|dx\leq \mathcal{C},  \quad \text{with } \quad \overline{v_p}^r =\fint_{\mathcal W_r(y_0)}v_p(y)dy.$$
Hereafter $\mathcal{C}$ denotes any constant depending only on $N$. Then, by \cite{JN}, there exists $\alpha>0$ independent of $N$ such that
$$\bigg(\fint_{\mathcal W_{\frac32}}e^{\alpha v_p}\bigg)\bigg(\fint_{\mathcal W_{\frac32}}e^{-\alpha v_p}\bigg)\leq \bigg(\fint_{\mathcal W_{\frac32}} e^{\alpha |v_p- \overline{v_p}^{3/2}|} \bigg)^2\leq \mathcal{C}.$$
That is,
\begin{equation}\label{Integral_comparison_32}
I\left((p-1)\alpha,3/2\right)\leq \mathcal{C}^{\frac{1}{p-1}}I\left((1-p)\alpha,3/2\right).
\end{equation}

Now, by the Sobolev inequality, for any $s\in\R$ and non-negative function $\varphi\in C_0^1(\mathcal W_2)$, we have
\begin{align}\label{sobolev-ineq-yield}
\left(\int_{\mathcal W_2}\varphi^{p^\ast}u_p^{sp^\ast}\right)^{\frac{1}{p^\ast}}&\leq \mathcal{C}\left(\int_{\mathcal W_2}F^p(\nabla(\varphi u_p^s))\right)^{\frac{1}{p}}\\
& \leq \mathcal{C}\left(\int_{\mathcal W_2}F^p(\nabla \varphi) u_p^{sp}\right)^{\frac{1}{p}}+\mathcal{C}|s|\left(\int_{\mathcal W_2}\varphi^pu_p^{sp-p}F^p(\nabla u_p)\right)^{\frac{1}{p}},
\end{align}
where $p^\ast=\frac{Np}{N-p}$. Moreover, if $s<0$, and as by \eqref{zVsz} $-\z \in \partial F(\nabla u_p)$, then 
\begin{align*}
\int_{\mathcal W_2}\varphi^pu_p^{sp-p}F^p(\nabla u_p)&=\int_{\mathcal W_2}\varphi^pu_p^{sp-p}F^{p-1}(\nabla u_p)(-\z)\cdot \nabla u_p\\
& =\frac{1}{sp-p+1}\int_{\mathcal W_2}\varphi^p F^{p-1}(\nabla u_p)(-\z)\cdot \nabla u_p^{sp-p+1}\\
&= \frac{p}{|sp-p+1|}\int_{\mathcal W_2}\varphi^{p-1}u_p^{sp-p+1}F^{p-1}(\nabla u_p)(-\z)\cdot \nabla\varphi,
\end{align*}
which follows by taking $u_p^{s p - p + 1} \varphi^p$ as a test function in \eqref{weakU}. Now using  \eqref{CS-ineq} and $F^\circ(-\z)\leq 1$ it holds
\begin{align*}
\int_{\mathcal W_2}\!\!\varphi^pu_p^{sp-p}F^p(\nabla u_p) \leq \tfrac{p}{|sp-p+1|}\left(\!\int_{\mathcal W_2}\!\!\varphi^{p}u_p^{sp-p}F^p(\nabla u_p)\!\right)^{\frac{1}{p'}}\!\!\left(\!\int_{\mathcal W_2}\!\!F^p(\nabla \varphi)u_p^{sp}\!\right)^{\frac1p}.
\end{align*}
Therefore,
$$\int_{\mathcal W_2}\varphi^pu_p^{sp-p}F^p(\nabla u_p)\leq \frac{p^p}{|sp-p+1|^p}\int_{\mathcal W_2}F^p(\nabla \varphi)u_p^{sp},$$
which, together with \eqref{sobolev-ineq-yield}, gives
$$\left(\int_{\mathcal W_2}\varphi^{p^\ast}u_p^{sp^\ast}\right)^{\frac{1}{p^\ast}}\leq \mathcal{C} \left(1+\frac{|s|p}{|sp-p+1|}\right)\left(\int_{\mathcal W_2}F^p(\nabla \varphi)u_p^{sp}\right)^{\frac1p}.$$
Now, for $1\leq r<R\leq 2$, choosing $\varphi$ appropriately we get
$$\left(\int_{\mathcal W_r}u_p^{sp^\ast}\right)^{\frac{1}{p^\ast}}\leq \frac{\mathcal{C}}{R-r} \left(1+\frac{|s|p}{|sp-p+1|}\right)\left(\int_{\mathcal W_R}u_p^{sp}\right)^{\frac1p}.$$
Hence, letting $\theta:=N/(N-p)$ and noticing that $\frac{|s|p}{|s p-p+1|}\leq 1$ for $s < 0$,
\begin{equation}\label{recursive_start}
  I(\theta sp,r)^s= I(sp^\ast,r)^s\leq \frac{ 2 \, \mathcal{C}}{R-r} I(sp,R)^s.
\end{equation}

Set $s_0:=(\alpha(1-p))/p$ and $r_0:=3/2$. For $k\in\N$, let $s_k:=\theta^ks_0$ and $r_k=1+2^{-(k+1)}$. Letting $s=s_k$, $r=r_{k+1}$ and $R=r_k$ in \eqref{recursive_start}, since $s_{k+1}=\theta s_k$, $r_k - r_{k+1}=2^{-(k+2)}$ and $s_k<0$, we obtain
\[I(s_{k+1}p,r_{k+1})\geq \mathcal{C}^\frac{1}{s_k}2^\frac{k+3}{s_k} I(s_kp,r_k)\]
for every $k$. Consequently,
\begin{align}
I((1-p)\alpha,3/2)&=I(s_0p,r_0) \leq 2^{-\frac{3}{s_0}}\mathcal{C}^{-\frac{1}{s_0}} I(s_{0}p,r_{1})
\\ & \leq 2^{-\frac{3}{s_0}-\frac{4}{s_1}}\mathcal{C}^{-\frac{1}{s_0}-\frac{1}{s_1}}I( s_{2}p,r_{2})\leq \cdots \leq 2^\gamma \mathcal{C}^\beta I(-\infty,1),
\end{align}
where
$$\beta=-\sum_{k=0}^\infty\frac{1}{s_k}=\frac{N}{\alpha(p-1)}$$
and, as $p < N$,
$$\gamma=-\sum_{k=0}^{\infty}\frac{k+3}{s_k}=\frac{N^2+2Np}{\alpha p(p-1)} \leq \frac{3 N}{\alpha (p-1)}  .$$
That is,
$$I((1-p)\alpha,3/2)\leq \mathcal{C}^{\frac{1}{p-1}}I(-\infty,1),$$
which, together with \eqref{Integral_comparison_32}, gives
$$I((p-1)\alpha,3/2)\leq e^{\frac{\mathcal{C}}{p-1}}I(-\infty,1).$$
This is equivalent to
$$\sup_{\mathcal W_1} v_p\leq \mathcal{C}-\frac{1}{\alpha}\log\bigg(\fint_{\mathcal W_{\frac{3}{2}}}e^{-\alpha v_p}\bigg).$$

We want to consider positive values of $s$ too. However, the factor
$$g(s):=\frac{|s|p}{|sp-p+1|}$$
appearing in the previous computations is unbounded near $\frac1{p'}=\frac{p-1}{p}$. Therefore, set
$$\sigma_0=\frac{2N-p}{2Np'}$$
so that
$$\sigma_0<\frac1{p'}<\sigma_1=\theta \sigma_0=\frac{2N-p}{2(N-p)p'}$$
and $g(\sigma_0)=g(\sigma_1)=2N/p-1$. Hence, for any $l\in\Z$ and $s=\theta^l \sigma_0$ we have
$$g(s)\leq \frac{2N}{p}-1.$$
Let $k_0\in\Z$ such that
$$\theta^{k_0}\sigma_0\leq \frac{\alpha}{p'}<\theta^{k_0+1}\sigma_0.$$
For $s_k:=\theta^{k_0+k}\sigma_0$ and $r_k:=1+2^{-(k+1)}$, we may repeat the arguments above to get
$$\sup_{\mathcal W_1} u_p= I(\infty,1)\leq e^{\frac{\mathcal{C}}{p-1}}I(-\infty,1) = e^{\frac{\mathcal{C}}{p-1}} \inf_{\mathcal W_1} u_p.$$
This is equivalent to
$${\rm osc}_{\mathcal W_1} v_p \leq \mathcal{C}.$$
Since, by the proof of Theorem \ref{cont-thm}, $v_p$ converges (up to a subsequence, not relabelled) locally uniformly to $v$, we obtain that $${\rm osc}_{\mathcal W_1} v \leq \mathcal{C}.$$\end{proof}

\section{About the anisotropic perimeter of sublevel sets} \label{sublevel-1}

 Concerning notation, all the sets appearing in this and the subsequent sections are considered to be Borel measurable subsets of $\R^N$.

\subsection{Anisotropic perimeter, densities and coarea formula}

 A set  $E$ within an open subset $U \subset \R^N$ is said to be of finite perimeter in $U$ if $\chi_{_E} \in BV(U)$. In this case, $P_{_{\! F}}(E;U)$ denotes the anisotropic perimeter (or $F$-perimeter, for brevity) of $E$ in  $U$, and is defined as
\begin{equation} \label{Fper-def}
	P_{_{\! F}}(E;U) :=  \big|D \chi_{_E}\big|_F(U).
\end{equation}
If $U = \R^N$ we simply write $P_{_{\! F}}(E) = P_{_{\! F}}(E;\R^N)$.  We say that $E$ has finite $F$-perimeter in $U$ if $P_{_{\! F}}(E;U) < \infty$. In turn, $E$ has locally finite $F$-perimeter in $U$ if, for every open set $B \Subset U$, it holds that $P_{_F}(E; B) < \infty$. By \eqref{FequivEuclid}, the finiteness of Euclidean and anisotropic perimeter are equivalent.

If $E$ has locally finite perimeter, it makes sense to consider the reduced  boundary
\[\partial^\ast E \!=\Big\{x \in \R^N \! :\,  \lim_{r \searrow 0} \frac{D \chi_{_E}(B_r(x))}{\big|D \chi_{_E}\big|(B_r(x))}\!=: \nu_{_E}(x) \text{ exists and belongs to } \mathbb{S}^{N-1}\Big\},\]
where $\nu^E(x)$ is known as  the measure-theoretic normal to $E$ at $x$, and clearly $\nu_{_E} = \nu_{\chi_{_E}}$. Then by de Giorgi structure theorem (see e.g. \cite[Theorem 5.9]{Maggi}) and \eqref{anis-decomp}, it is known that
\begin{equation} \label{Fper}
P_{_{\! F}}(E;U) = \int_{U\cap\partial^\ast E} F(\nu^E) \, d \mathcal H^{N-1}.
\end{equation}

Recall that, for $\alpha \in [0,1]$,  the set of points of density $\alpha$ in $E$ is given by
\[E^{\alpha}:=\bigg\{x \in \R^N \, : \, \lim_{r \searrow 0} \frac{\big|B_r(x) \cap E\big|}{|B_r(x)|} = \alpha\bigg\}.\]
Typically, $E^1$ is known as measure theoretic interior of $E$. Hereafter $E\Delta A := (E \setminus A) \cup (A \setminus E)$ denotes the symmetric difference of the sets $E$ and $A$.

\begin{lemma} \label{prop-int-0} If $E$ is a set of finite perimeter, the following properties hold:
	{\rm (a)} {\bf (Federer)} \cite[Theorem 3.61]{Ambrosio}  $\partial^\ast E \subset E^{1/2}$. Moreover, $E$ has density either 0 or 1/2 or 1  $\mathcal H^{N-1}$-a.e. in $\R^N$.
	
	{\rm (b)} \cite[Example 5.17 and Remark 15.2]{Maggi}  $\big|E\Delta E^{1}\big| = 0$.  Consequently, $P_{_F}(E) = P_{_F}(E^{1})$ and $\partial^\ast E = \partial^\ast E^1$.
\end{lemma}

 The goal of this section is to deal with the sublevel sets $E_t, G_t \subset \Omega^e$ of solutions of A-IMCF, defined as in \eqref{def-sublevel}, and show that their $F$-perimeters are finite and coincide.

Notice that we can extend the solution $v$ of \eqref{weak-IMCF} to be identically zero within $\overline{\Omega}$ (in fact, the boundary condition in \eqref{IMCF-bdry} guarantees that this is well-defined and hence the extended function belongs to $DBV_{\rm loc}(\R^N) \cap L^\infty_{\rm loc}(\R^N)$). We will use the extension when needed without further mention, but the corresponding vector field $\z$ will not be extended.

 Accordingly, $G_t$ and $E_t$ are redefined as
\[G_t =\{x \, :\, \bar v(x) \leq t\}  \quad \text{for} \quad \bar v(x) = \left\{\begin{array}{ll} v(x), & \text{if } x \in \Omega^e, \smallskip \\ 0, &\text{if } x \in \Omega,\end{array}\right.\]
so that $\overline\Omega \subset G_t$ for all $t \geq 0$, and similarly $E_t =\{x \, :\, \bar v(x) < t\}$.

\begin{remark} \label{Dv=0-on-v=t}
Observe that, as the measure $|Dv|_{_F}\res\{v=t\}$ vanishes and $v$ has no jump, by the previous extension, it also holds
\begin{equation} \label{TV-G0-0}
|Dv|_{_F}(\{v=0\})= |Dv|_{_F}(G_0) = 0.
\end{equation}
\end{remark}

As the solution $v$ from Theorem \ref{exist-thm} belongs to $BV_{\rm loc}({\R^N})$,  the anisotropic version of the coarea formula \cite[Remark 4.4]{AmBe} yields 
\begin{align} \label{F-coarea}
	|Dv|_{_{\! F}}(U) = \int_\R P_{_{\! F}}(\{v > t\}; U)\, dt \quad \text{for every open set } U\subset\R^N.
\end{align}
 In particular, this ensures that the measures $\!\displaystyle \int_\R \!\big|D \chi_{_{\{v > t\}}}\big|_F(A)\, dt$ and $|Dv|_{_F}(A)$ agree for any open subset $A$ of $U$, and thus coincide. Hence, for any $B \subset U$,
\begin{equation} \label{coarea-2}
|Dv|_{_F}(B) = \int_\R \big|D \chi_{_{\{v > t\}}}\big|_F(B)\, dt.
\end{equation}

Notice that \eqref{F-coarea} guarantees that the superlevel sets $\{v > t\}$, and thus $G_t$ have locally finite $F$-perimeter for a.e. $t \geq 0$. On the other hand, since $v$ tends to $\infty$ as $|x|\to \infty$, we can find $r(t) > 0$ such that
\begin{equation} \label{GsubB}
G_t \subset B_{r(t)}(0)
\end{equation}
and hence the sublevel sets $G_t$ have finite $F$-perimeter for a.e. $t >0$. The goal is to show that they actually have finite $F$-perimeter for all $t >0$.

\subsection{Reduced boundary of the sublevel sets}

We need an auxiliary result about the reduced bounda\-ry of the sublevel sets. Let us first prove it for a.e. $t \geq 0$.
\begin{lemma} \label{v=t}
	Let $t \geq 0$ be such that $G_t$ has finite $F$-perimeter. Then
	 \[v = t \quad  \ \mathcal H^{N-1}\text{-a.e. in }  \partial^\ast G_t.\]
Therefore, given $0\leq \tau<t$ such that $G_\tau$ and $G_t$ have finite $F$-perimeter, $\mathcal{H}^{N-1}(\partial^\ast G_t\cap\partial^\ast G_\tau)=0$; in particular, $\mathcal{H}^{N-1}(\partial^\ast G_t\cap\partial^\ast \Omega)=0$, leading to
\begin{equation} \label{eq-per}
P_{_F}(G_t; \Omega^e) = P_{_F}(G_t).
\end{equation}
	
\end{lemma}

\begin{proof}
	Since $\mathcal H^{N-1}(S_v) = 0$ as there is no jump part (cf. Lemma \ref{no-jump}), and by Lemma \ref{prop-int-0} (a), it is enough to prove the claims for Lebesgue points of density $1/2$, that is, for $x$ so that
	\begin{equation} \label{density12}
	\lim_{r \searrow 0} \frac{\big|B_r(x) \cap G_t\big|}{|B_r(x)|} = \frac1{2} \quad \text{and} \quad \lim_{r \searrow 0} \frac{\big|B_r(x) \cap \{v>t\}\big|}{|B_r(x)|} = \frac1{2}.
	\end{equation}
	The first equality leads to
	\[\lim_{r \searrow 0}\fint_{B_r(x) \cap G_t} |v - v(x)| \leq 2\lim_{r \searrow 0}\fint_{B_r(x) } |v - v(x)| =0,\]
	where the latter follows because $x$ is a Lebesgue point of $v$ (recall Definition \ref{lebesgue}). In particular,
	\[v(x) =\lim_{r \searrow 0}\fint_{B_r(x) \cap G_t} v \leq t.\]
	By an analogous argument but with the second equality in \eqref{density12}, we reach
	\[v(x) = \lim_{r \searrow 0}\fint_{B_r(x) \cap \{v > t\}} v \geq t,\]
	which finishes the proof of the first claim. Now \eqref{eq-per} follows directly from the fact that $\partial G_t \cap \overline \Omega = \emptyset$.
\end{proof}

\subsection{Anisotropic perimeter of the sublevel sets}

\begin{lemma} \label{lem:sublevel}
Let $v$ be a weak solution in the sense of Definition \ref{weak-IMCF-def}. Then the sublevel sets $G_t$  have finite anisotropic perimeter for every $t \geq 0$.
\end{lemma}
\begin{proof}
	By the general coarea formula in \cite[Theorem 4.2]{CraCi} and equation \eqref{weak-IMCF}, we can write, for any  set $B \subset \Omega^e$,
	\begin{align}
	\int_{\R} (\z, D \chi_{_{\{v > t\}}})(B) \, dt  &= (\z, Dv) (B) = |Dv|_{_{\! F}}(B) 
	 = \int_{\R}  \big|D \chi_{_{\{v > t\}}}\big|_F(B) \, dt,
	\end{align}
 where the latter follows by \eqref{coarea-2}. From here, since $-(\z, D \chi_{_{G_t}}) \leq \big|D \chi_{_{G_t}}\big|_{_F}$, we have the equality
\begin{equation} \label{meas-eq1}
-(\z, D \chi_{_{G_t}}) = \big|D \chi_{_{G_t}}\big|_{_F} \quad \text{as measures in } \Omega^e \ \text{ for a.e. } \ t > 0.
\end{equation}	

	{\bf Claim 1.}  $P_{_{\! F}}(G_t) - |Dv|_{_{\! F}}(G_t) = \Lambda$, where $\Lambda$ is constant, for a.e. $t > 0$.

Take two times $T > \tau > 0$ satisfying \eqref{meas-eq1} and so that $G_T$ and $G_\tau$ have finite $F$-perimeter,  and hence  \eqref{eq-per} holds. 
Now by substituting \eqref{meas-eq1} in \eqref{Fper-def} and applying Green's formula in $B_{r(T) +1}(0) \cap \Omega^e$, where $r(T)$ is the radius from \eqref{GsubB}, we get
\begin{align*}
P_{_{\! F}}(G_T) - P_{_{\! F}}(G_\tau) & = -\int_{\Omega^e} (\z, D \chi_{_{G_T}}) + \int_{\Omega^e} (\z, D \chi_{_{G_\tau}}) \\ & = -\int_{\Omega^e \cap B_{r(T)+1}(0)} (\z, D \chi_{_{\{\tau < v \leq T\}}}) \\ & = \int_{\Omega^e \cap B_{r(T)+1}(0)} \chi^\ast_{_{\{\tau < v \leq T\}}} d (\divi\z) = \int_{{\{\tau < v \leq T\}}} d(\divi \z),
\end{align*}
where for the last equality we use that $\chi^\ast_{_{\{\tau < v < T\}}} = 1$  on $\{\tau<v<T\}$ and Corollary \ref{div=0}. Hence equation \eqref{weak-IMCF} yields
\begin{align*}
|Dv|_{_{\! F}}(G_T)- |Dv|_{_{\! F}}(G_\tau) & = \int_{\{\tau < v \leq T\}} d |Dv|_{_{\! F}} =  \int_{{\{\tau < v \leq T\}}} d(\divi \z) \\ & = P_{_{\! F}}(G_T) - P_{_{\! F}}(G_\tau).
\end{align*}
Thus the claim follows for a.e. $t >0$ because we argued for any pair of times satisfying \eqref{meas-eq1}.

{\bf Claim 2.}  $\displaystyle \frac{d^+}{dt} |Dv|_{_{\! F}}(G_t) = \Lambda + |Dv|_{_{\! F}}(G_t)$ for every $t \geq 0$, where $\frac{d^+}{dt}$ denotes the right derivative.

First, notice that by  \eqref{coarea-2} we have
\[|Dv|_{_{\! F}}(G_t) = \! \int_0^t \!\!\int_{\R^N}   \big|D \chi_{\{v > s\}}\big|_{_F} \, ds = \! \int_0^t \!\! \int_{\R^N}  \big|D \chi_{\{v \leq  s\}}\big|_{_F} \, ds = \! \int_0^t \!\! P_{_{\! F}}(G_s) \, ds.\]
In particular,
\begin{equation} \label{abs-cont}
t \mapsto |Dv|_{_{\! F}}(G_t) \ \text{ is absolutely continuous on each bounded interval. }
\end{equation}
Moreover, for any $h >0$  and any $t \geq 0$, Claim 1 leads to
\begin{align*}
\frac{|Dv|_{_{\! F}}(G_{t+h})- |Dv|_{_{\! F}}(G_t)}{h} = \frac1{h} \int_t^{t+h} P_{_{\! F}}(G_s) \, ds = \Lambda +  \frac1{h} \int_t^{t+h} |Dv|_{_{\! F}}(G_{s}) \, ds.
\end{align*}
Taking limits as $h \to 0^+$, this implies that there exists the right derivative and the claimed equality follows.

	{\bf Claim 3.}  $P_{_{\! F}}(G_t) \leq \dfrac{d^+}{dt} |Dv|_{_{\! F}}(G_t)$ for every $t \geq 0$.

Fix $t \geq 0$ and $h > 0$, for which we consider the Lipschitz continuous function $\eta: \R \rightarrow [0,1]$ given by
\begin{equation} \label{def-eta}
\eta(s) = \left\{\begin{array}{ll} 1, & \text{if } s \leq t, \smallskip \\  1 + \dfrac{t-s}{h}, & \text{if } s \in [t, t+h], \smallskip \\ 0, & \text{if } s \geq t+h. \end{array} \right.
\end{equation}	
Notice that, for $s \neq t, t+h$, it holds $\eta'(s) = - \frac1{h} \chi_{(t, t+h)}(s)$.	Using this, together with Corollary \ref{div=0}, we reach
\begin{align*}
\frac{|Dv|_{_{\! F}}(G_{t+h})- |Dv|_{_{\! F}}(G_t)}{h} & = \frac1{h} \int_{\{t < v< t+h\}} dF(Dv) = \int_{\R^N} -\eta'(v) dF(Dv) \\ &= \int_{\R^N} dF(D(\eta \circ v)) \geq \int_{\R^N} \eta(v) \divi {\bf w},
\end{align*}
for any test function ${\bf w} \in C^1_{c}(\R^N; \R^N)$ with $F^\circ({\bf w}) \leq 1$. Notice that the first equality in the second line follows from the chain rule in Lemma \ref{DBV} (a), while the last inequality comes from \eqref{FtotV}.

Accordingly, taking limits as $h \to 0^+$, we have  that
\[\frac{d^+}{dt} |Dv|_{_{\! F}}(G_t) \geq \int_{\R^N} \chi_{_{G_t}} \divi {\bf w}   \]
and, taking the supremum over all such vector fields ${\bf w}$,
\[\frac{d^+}{dt} |Dv|_{_{\! F}}(G_t) \geq  \big|D \chi_{_{G_t}}\big|_{_F}(\R^N) = P_{_F}(G_t),\]
thus the claimed inequality follows.

Summing up, if we join Claims 2 and 3, by means of \eqref{GsubB} one can conclude
\[P_{_{\! F}}(G_t) \leq \Lambda + |Dv|_{_{\! F}}(G_t) < \infty \quad \text{for every } t \geq 0.\]
\end{proof}

\begin{remark} \label{rmk-interior}
Now the conclusions of Lemma \ref{v=t} are valid for all times. Moreover, as we are interested in properties involving the $F$-perimeter, which are invariant under modifications of 0-Lebesgue measure by Lemma \ref{prop-int-0}, hereafter we will identify any finite perimeter set (in particular, $G_t$ for any $t\geq 0$) with its measure theoretic interior.
\end{remark}

\smallskip

\begin{lemma} \label{lem:sublevel-2}
	Let $v$ be the weak solution from  Definition \ref{weak-IMCF-def}. Then the sublevel sets $E_t$ have finite anisotropic perimeter $P_{_{\! F}}$ for every  $t > 0$.
\end{lemma}

\begin{proof}
Arguing exactly as before, but for $h < 0$, one gets that, for the left derivative and the same constant from Claim 1 above, it holds	
\[\displaystyle \frac{d^-}{dt} |Dv|_{_{\! F}}(G_t) = \Lambda + |Dv|_{_{\! F}}(G_t).\]

One can also mimic the proof of Claim 3, but using the following function instead of the one introduced in \eqref{def-eta}:
\begin{equation}
\xi(s) = \left\{\begin{array}{ll} 1, & \text{if } s \leq t+h, \smallskip \\   \dfrac{s-t}{h}, & \text{if } s \in [t+h, t], \smallskip \\ 0, & \text{if } s \geq t. \end{array} \right.
\end{equation}	
Doing so, we conclude
\[P_{_{\! F}}(E_t) \leq \dfrac{d^-}{dt} |Dv|_{_{\! F}}(G_t) \quad \text{for all } t > 0,\]
and arguing as before the statement follows.
\end{proof}

\subsection{The generalized unit outward normal and continuity of $F$-perimeter}

Here we show that, in case of extra regularity, the vector field $\z$ corresponding to the weak solution $v$ of \eqref{weak-IMCF} represents the generalized outward unit normal to $\partial^\ast G_t$. Recall that
	\[\nu_{_{G_t}} = \nu_{_{\partial^\ast G_t}} = \nu_v \quad \mathcal{H}^{N-1}\text{-a.e. in} \ \partial^\ast G_t.\]

\begin{proposition} \label{meas-Gt}
	With the same notation of the previous subsections, we have the equality
	\begin{equation} \label{meas-char-eq}
	-(\z, D \chi_{_{G_t}}) = \big|D \chi_{_{G_t}}\big|_{_F} \quad \text{as measures on } \  \Omega^e\  \text{ for all } \  t \geq 0.
	\end{equation}
In particular, if  $\z \cdot \nu_v$ is continuous, this leads to
\[\z \cdot \nu_v = -F(\nu_v) \quad \mathcal H^{N-1}\text{-a.e. in } \partial^\ast G_t.\]
\end{proposition}

\smallskip

\begin{remark}
Notice that the last equality means that $-\z \in \partial F(\nu_v)$. In case that $F$ is smooth, this implies that $\z$ points in the direction of the anisotropic outward normal to $\partial^\ast G_t$ $\mathcal H^{N-1}$-a.e., and then the left hand side of \eqref{weak-IMCF} is actually the anisotropic mean curvature $H_F = \divi \z$. Hence our definition of A-IMCF recovers that of \cite{IAMCF}.
\end{remark}

\begin{proof}
Consider the function $\eta$ defined in \eqref{def-eta} for any $t \geq 0$. By means of Green's formula applied on $\Omega_{r}:= \Omega^e \cap B_{r(t)+1}$ so that \eqref{GsubB} holds, we get
\[\int_{\Omega_{r}} (\eta \circ v) \, d(\divi \z) = - \int_{\Omega_{r}} \big(\z, D(\eta \circ v)\big) + \int_{\partial \Omega} [z, \nu],\]
where we have applied that $\eta(v) = 1$ on $\partial \Omega$ since $v$ vanishes on the boundary.

Now, as $-\eta$ is a Lipschitz non-decreasing function, by Lemma \ref{DBV} (b1),
\[\big(\z, D(-\eta \circ v)\big) = -\eta'(v) (\z, Dv) \quad \text{as measures.}\]
Hence, by means of equation \eqref{weak-IMCF}, we obtain
\begin{align*}
\int_{\Omega_r} (\eta \circ v) \, d(\divi \z) - \int_{\partial \Omega} [\z, \nu] & = \int_{\Omega_r} \!\!\!\!-\eta'(v) dF(Dv) = \frac1{h} \int_{\{t < v <t+h\}}  \!\!dF(Dv)
\\ & =  \frac{|Dv|_{_{\! F}}(G_{t+h})- |Dv|_{_{\! F}}(G_t)}{h},
\end{align*}
where we applied Corollary \ref{div=0}. Letting $h \to 0^+$ yields by Green's formula
\begin{align} \label{aux-per}
\displaystyle \frac{d^+}{dt} |Dv|_{_{\! F}}(G_t) & = \int_{\Omega_r} \chi_{_{G_t}}^\ast d(\divi \z) - \int_{\partial \Omega} [\z, \nu] \nonumber \\ & = - \int_{\Omega^e} d(\z, D\chi_{_{G_t}}) \leq \int_{\Omega^e} dF(D \chi_{_{G_t}}) \leq  P_{_{\! F}}(G_t),
\end{align}
which follows from the Cauchy-Schwarz type inequality in Lemma \ref{prop-ATV-2} (a) because $F^\circ(\z) \leq 1$.
This means that we have equality in Claim 3 within the proof of Lemma \ref{lem:sublevel}, that is,
\begin{equation} \label{eq-PF-TV}
P_{_F}(G_t) = \dfrac{d^+}{dt} |Dv|_{_{\! F}}(G_t) = \Lambda + |Dv|_{_{\! F}}(G_t) \quad \text{for all } t\geq 0.
\end{equation}
 Hence we get equality also in \eqref{aux-per}, that is, $0=    \int_{\Omega^e}   d(-\z, D\chi_{_{G_t}}) - dF(D \chi_{_{G_t}})$ which, again by Cauchy-Schwarz from Lemma \ref{prop-ATV-2} (a), proves the first assertion of the statement.

For the second claim about the vector field $\z$,   when $\z \cdot \nu_v$ is continuous,
  by using \eqref{pairing-jump} and \eqref{anis-decomp}, we have
 \[\z \cdot  \nu_v = [\z, \nu_v]= \frac{[\z, \nu_v]^+ + [\z, \nu_v]^-}{2} = - F(\nu_v)  \quad \mathcal H^{N-1}\text{-a.e. in } \partial^\ast G_t,\]
 where the first equality comes from \cite[Corollary 4.6]{CoCraCiMa}.
\end{proof}

\subsection{Continuity and exponential growth of the F-perimeter}

 Next, we show that, despite possible {\it fattening} of the weak solution of the A-IMCF, the $F$-perimeter is preserved.

\begin{corollary}[Continuity of the anisotropic perimeter]
	Let $v$ be a weak solution of the A-IMCF in the sense of Definition \ref{weak-IMCF-def}, then the function $t \mapsto P_{_{\! F}}(G_t)$ is right-continuous for every $t \geq 0$ and $t \mapsto P_{_{\! F}}(E_t)$ is left-continuous for every $t > 0$. Moreover, it holds
	\begin{equation} \label{cont-per}
	P_{_{\! F}}(G_t) = P_{_{\! F}}(E_t) \qquad \text{for all} \ t >0.
	\end{equation}
\end{corollary}

\begin{proof}
	First, notice that by \eqref{eq-PF-TV}, we have that for any $s \geq 0$
\[P_{_F}(G_s) = \Lambda + |D v|_{_F}(G_s),\]
and taking limits as $s\searrow t$, thanks to \eqref{abs-cont}, we  get
\[\lim_{s \searrow t} P_{_F}(G_s) = \Lambda + |D v|_{_F}(G_t) = P_F(G_t), \]
which gives the claimed right-continuity.

Repeating the arguments from the proof of Proposition \ref{meas-Gt} for the sublevel sets $E_t$, one can show that
\begin{equation} \label{eq-PF-TV-2}
P_{_F}(E_t) = \dfrac{d^-}{dt} |Dv|_{_{\! F}}(G_t) = \Lambda + |Dv|_{_{\! F}}(G_t) \quad \text{for all } t > 0,
\end{equation}
which leads to the claimed equality thanks to  \eqref{eq-PF-TV}. The left continuity follows arguing as before.
\end{proof}

 Summing up, if we introduce the functional
	\begin{equation}
	\mathcal J_v(E):=P_{_{\! F}}(E) - |Dv|_{_F}(E) \quad \text{for a bounded set } E\subset \R^N,
	\end{equation}
we conclude from \eqref{eq-PF-TV} that
\begin{equation} \label{JE=JG}
P_{_F}(G_0) = \mathcal J_v(G_0) = \Lambda = \mathcal J_v(G_t)=\mathcal J_v(E_t) \qquad  \text{for every } \ t >0.
\end{equation}
Here the first equality comes from \eqref{TV-G0-0}, and the last one using Corollary \ref{div=0}, as well as \eqref{cont-per}.

As in the classical (smooth and isotropic) case, it follows that the rescaled $F$-perimeter $e^{-t} P_{_F}(G_t)$ stays constant.

\begin{corollary}
Let $v$ be a weak solution of the A-IMCF in the sense of Definition \ref{weak-IMCF-def}, then it holds
\begin{equation} \label{exp-per}
P_{_{\! F}}(G_t) =  e^tP_{_{\! F}}(G_0)  = e^t P_{_{\! F}}(\{v=0\}) \quad \text{for all } \ t > 0.
\end{equation}

\end{corollary}

\begin{proof}

By \eqref{div-exp0}, we have that ${\divi (e^{-v}\z)} = 0$ as measures in $\Omega^e$, and hence by application of Green's formula in $B_{r(t)+1}(0) \cap \Omega^e$ satisfying \eqref{GsubB}, we reach
\[0  = - \int_{\{s < v \leq t\}} \divi(e^{-v} \z) = \int_{\Omega^e} \big(e^{-v} \z, D \chi_{_{\{s < v \leq t\}}}\big) = \int_{\Omega^e} e^{-v}\big(\z, D \chi_{_{\{s < v \leq t\}}}\big).\]
where the latter follows from \eqref{prule-pair}. As $\chi_{_{\{s < v \leq t\}}} = \chi_{_{G_t}} - \chi_{_{G_s}}$, \eqref{meas-char-eq} yields
\begin{align} \label{aux-eq7.6}
0 & = \int_{\Omega^e} e^{-v} \big(\z, D\chi_{_{G_t}}\big) -  \int_{\Omega^e} e^{-v} \big(\z, D\chi_{_{G_s}}\big)\nonumber
\\ & = - \int_{\Omega^e} e^{-v} \big|D\chi_{_{G_t}}\big|_{_F} +  \int_{\Omega^e} e^{-v} \big|D\chi_{_{G_s}}\big|_{_F}.
\end{align}

On the other hand, for $0 < s<t$, we can apply Lemma \ref{v=t} recalling that by \eqref{Fper} the Radon measure $\big|D\chi_{_{G_s}}\big|_{_F}$ has support in $\partial^\ast G_s$, and hence
\begin{align*}
e^{-s} P_{_{\! F}}(G_s) - e^{-t} P_{_{\! F}}(G_t) & = e^{-s} \int_{\R^N} \big|D\chi_{_{G_s}}\big|_{_F}-  e^{-t} \int_{\R^N} \big|D\chi_{_{G_t}}\big|_{_F} \\ & =  \int_{\R^N} e^{-v}\big|D\chi_{_{G_s}}\big|_{_F}-  \int_{\R^N} e^{-v} \big|D\chi_{_{G_t}}\big|_{_F} = 0,
\end{align*}
which comes from \eqref{aux-eq7.6}. The statement follows by letting $s \to 0^+$.
\end{proof}

\section{Outward minimizing properties of solutions} \label{out-min}

\subsection{Variational definition of weak solutions}

As in \cite{HuIl, Moser15} we consider a variational version of the definition of weak solutions, and we will show that the solution from Definition \ref{weak-IMCF-def} is also a solution in the following sense:

\begin{definition}
	A function $v \in BV_{\rm loc}(\Omega^e) \cap L^\infty_{\rm loc}(\Omega^e)$ is a weak variational solution of A-IMCF if for every compact set $K \subset \Omega^e$ and for every  $w \in BV_{\rm loc}(\Omega^e) \cap L^\infty_{\rm loc}(\Omega^e)$ with $w = v$ in $\Omega^e \setminus K$ we have
	\[\mathcal I_v^K(v) \leq  \mathcal I^K_v(w),\]
	where $\mathcal I_v^K$ is the functional given by
	\[\mathcal I_v^K(w):= \int_K dF(Dw) + \int_K w \, dF(Dv).\]
\end{definition}

\begin{proposition} \label{def-var}
	The weak solution $v$ in the sense of Definition \ref{weak-IMCF-def} is a weak variational solution of A-IMCF.
\end{proposition}

\begin{proof}
	Fix any function $w \in BV_{\rm loc}(\Omega^e) \cap L^\infty_{\rm loc}(\Omega^e)$ such that $\{w \neq v\}$ is contained in a compact set $K\subset\Omega^e$. We also take an open and bounded set $B$ with Lipschitz boundary and so that $K \subset B \subset \Omega^e$.
	
	We know that there exists a vector field ${\rm {\bf z}} \in \mathcal D \mathfrak M^\infty_{\rm loc}(\Omega^e)$ with $F^\circ({\rm {\bf z}}) \leq 1$ and satisfying ${\rm div} \, \z = F(Dv)$. If we multiply this equation by $w$, integrate over $B$ and apply Green's formula \eqref{green}, we get
	\begin{equation} \label{aux-var}
		-\int_B d(\z, Dw) + \int_{\partial B} [\z, \nu] w \, d \mathcal H^{N-1} = \int_B w \, dF(Dv).
	\end{equation}
	
	Now, Lemma \ref{prop-ATV-2} (a) ensures that $(z, Dw) \leq F(Dw)$, and hence
	\begin{align*}
		\int_{\partial B} [\z, \nu] w \, d \mathcal H^{N-1} \leq \int_B dF(Dw)+ \int_B w \, dF(Dv).
	\end{align*}

	On the other hand, arguing as before with $v$ instead of $w$, we have that
	\[\int_{\partial B} [\z, \nu] v \, d \mathcal H^{N-1} = \int_B d F(Dv)  + \int_B v \, dF(Dv),\]
	where we have also applied the second equality in \eqref{weak-IMCF}.
	
	Next, notice that $v = w$ on $\partial B$, which leads to
	\[\int_B d F(Dv)  + \int_B v \, dF(Dv) \leq \int_B dF(Dw)+ \int_B w \, dF(Dv).\]
	As this holds for all $B$, we have that $\mathcal I_v^K(v) \leq \mathcal I_v^K(w)$, as desired.
\end{proof}

\subsection{Geometric definition of weak solutions}
To achieve a geometric notion of weak solutions as in \cite{HuIl}, we
consider the minimization problem
\begin{equation}\label{minimize-geom}
  \inf \left\{ \mathcal J_v(E) \, : \, \Omega\subset E\Subset \R^N\right\}.
\end{equation}
The goal is to prove  the following result:
\begin{proposition} \label{equiv-sol}
	Let   $v \in BV_{\rm loc}(\Omega^e) \cap L^\infty_{\rm loc}(\Omega^e)$ be a weak solution of \eqref{EquationVp1}. Then, for each $\tau\geq 0$ and $t > 0$, $E_t$ and $G_\tau$  are minimizers of \eqref{minimize-geom}.
\end{proposition}

\begin{proof}
We start by proving that $G_t$ with $t> 0$ is a minimizer of
\begin{equation}\label{minimize-geom-overline}
	\inf \left\{ \mathcal J_v(E) \, : \, \overline\Omega\subset E\Subset \R^N\right\}.
\end{equation}
Let $E$ be a bounded set containing $\overline{\Omega}$. Then,
\begin{align} \label{diff-J}
\hspace*{-0.8cm} \mathcal J_v(E)-\mathcal J_v(G_t) & = P_{_{\! F}}(E)-P_{_{\! F}}(G_t) - |Dv|_{_F}(E\setminus G_t) + |Dv|_{_F}(G_t\setminus E).
\end{align}
In turn,  Lemma \ref{v=t}  ensures that $P_{_F}(G_t) = P_{_F}(G_t; \Omega^e)$, which implies by means of Cauchy-Schwarz from  Lemma \ref{prop-ATV-2} (b) that
\begin{align*}
	P_{_{\! F}}(E)-P_{_{\! F}}(G_t)
	& 
	 \geq  \big(\big|D \chi_{_E}\big|_{_F} - \big|D \chi_{_{G_t}}\big|_{_F}\big)(\Omega^e) \\
	& \geq -\int_{\Omega^e} d(\z, D \chi_{_E})+\int_{\Omega^e} d(\z, D \chi_{_{G_t}})
	\\
& =  -\int_{\Omega^e} d(\z, D \chi_{_E\setminus G_t})+\int_{\Omega^e} d(\z, D \chi_{_{G_t\setminus E}}).
\end{align*}
Now,  Green's formula on an open bounded set  with Lipschitz continuous boundary $B \supset G_t \cup E$ and  \eqref{weak-IMCF} allow us to conclude
\begin{align*}
P_{_{\! F}}(E)-P_{_{\! F}}(G_t)
	& = \int_{E\setminus G_t} d(\divi\z)-\int_{G_t\setminus E} d(\divi\z) \\
	& =  |Dv|_{_F}(E\setminus G_t)-|Dv|_{_F}(G_t\setminus E),
\end{align*}
which, by substitution into \eqref{diff-J}, leads to
\[\mathcal J_v(E)-\mathcal J_v(G_t)  \geq 0.\]

Let us now prove that, in fact, $G_t$ is a minimizer of \eqref{minimize-geom}. Given $E\supset\Omega$, consider $\tilde E=E\cup\partial\Omega\supset\overline\Omega$. Then, since $G_t$ is a minimizer of \eqref{minimize-geom-overline},
$$\mathcal J_v(G_t) \leq \mathcal J_v(\tilde E)=P_{_{\! F}}(E\cup\partial\Omega) - |Dv|_{_F}(E\cup\partial\Omega)=P_{_{\! F}}(E) - |Dv|_{_F}(E)= \mathcal J_v(E).$$
because $|\partial\Omega|=0$ and, by \eqref{TV-G0-0}, $|Dv|_{_F}(\partial\Omega)=0$.

Finally, the statement for $G_0$ and $E_t$ follows by \eqref{JE=JG}.
\end{proof}

\subsection{Definition and proof of the outward minimizing properties}

We first define the anisotropic analogue of the notion of minimizing hull or outward minimizing set.

\begin{definition}
	A measurable subset $E \subset \R^N$ of finite $F$-perimeter is said to be outward $F$-minimizing  if
	\[P_{_F}(E) \leq P_{_F}(A) \quad \text{for all} \quad  E \subset A  \Subset \R^N.\]
We say that it is strictly outward $F$-minimizing if the above inequality is strict unless $|A \setminus E| = 0$.
\end{definition}

 Notice that for smooth $F$ and $\partial E$ of class $C^2$, the above definition implies that $E$ is $F$-mean convex (i.e., $H_F \geq 0$).

 On the other hand, as pointed out in \cite[p.625]{Cha-Nov}, it holds

\begin{lemma} \label{prop-int}
	If $E$ is outward $F$-minimizing, then $E^1$ is open.
	
\end{lemma}
\noindent Recall that we consider $E^1$ as a representative of each finite perimeter set (cf.\ Remark \ref{rmk-interior}), thus we may assume that outward $F$-minimizing sets are open.

\begin{proposition} \label{min-hull}
	Let $v \in BV_{\rm loc}(\Omega^e) \cap L^\infty_{\rm loc}(\Omega^e)$ be a weak solution of the A-IMCF. Then
	\begin{enumerate}
		\item[{\rm (a)}] $E_t$ is outward $F$-minimizing for any $t > 0$.
		\item[{\rm (b)}] $G_t$ is strictly outward $F$-minimizing  for any $t \geq 0$.
	\end{enumerate}
\end{proposition}

\begin{proof}
	(a) Fix $t \geq 0$ and take any bounded $A \supset E_t$. By Proposition \ref{equiv-sol}
	\[P_{_F}(E_t) -  |Dv|_{_F}(E_t) \leq P_{_F}(A) - |Dv|_{_F}(E_t).\]
Hence
	\[P_{_F}(E_t)\leq  P_{_F}(E_t) + |Dv|_{_F}(A \setminus E_t)\leq P_{_F}(A)\]
	as desired.	Similarly, $G_t$ are outward $F$-minimizing for any $t \geq 0$.
	
	(b) Let us see that $G_t$ is strictly outward $F$-minimizing. With this aim, take a bounded set $A$ such that  $G_t \subset A$ and
	\[P_{_F}(G_t) = P_{_F}(A),\]
	 which implies that $A$ is also outward $F$- minimizing. Recall then that we take open representatives for both  $G_t$ and $A$. Arguing as before,
	\[P_{_F}(G_t)\leq  P_{_F}(G_t) + |Dv|_{_F}(A \setminus G_t)\leq P_{_F}(A).\]
	 This implies that $|Dv|_{_F}(A \setminus G_t) = 0$, which yields that $v$ is equivalent to a constant on each connected component of the open set $A \setminus G_t$.
	
	However, as $G_t$ is outward $F$-minimizing, the reduced boundary of these components intersect $\partial^*G_t$ on a set of positive $\mathcal H^{N-1}$-measure. Indeed, suppose otherwise that $C$ is a connected component of $A \setminus G_t$ with  $\mathcal H^{N-1}(\partial^* C\cap \partial^\ast G_t)=0$. Then, this yields $\mathcal H^{N-1}(\partial^*C\cap\partial^*(A\setminus C))=0$  which implies that
$$P_F(A\setminus C)=P_F(A)-P_F(C)<P_F(A)=P_F(G_t);$$
but $G_t \subset A\setminus C$, this is in contradiction with $G_t$ being outward $F$-minimizing.
	
	 Thus  $\mathcal H^{N-1}(\partial^* C\cap \partial^* G_t)>0$ for each connected component $C$ of $A\setminus\overline{G_t}$. Therefore, as $v \in DBV$ is  constant on each of these connected components and,  by Lemma \ref{v=t}, $v=t$ $\mathcal H^{N-1}$-a.e.\ on $\partial^\ast G_t$, we conclude that $v = t$ on $A \setminus G_t$. Consequently, $A \subset G_t$  and, as we are assuming the other inclusion, equality holds, as desired.
\end{proof}

Now we introduce the notion of $F$-minimal hull  (strictly outward $F$-minimizing hull) of a bounded set $E\subset\R^N$.
\[E' = \bigcap\big\{S \, : \, S \text{ bounded, strictly outward $F$-minimizing and } S \supset E\big\}.\]

 \begin{remark}
As we are working up to measure zero modification, the above definition may be realized by a countable intersection. Then arguing as in \cite[Proposition 1.3]{Tama}, as the proof basically relies on subadditivity of perimeters, one gets that	 $E'$ is strictly outward $F$-minimizing.
	\end{remark}

\begin{corollary}
With the same assumptions from Proposition \ref{min-hull}, it holds
\[\Omega'=G_0 \ \text{ and } \ E_t' = G_t \quad \text{up to null sets and for all } t > 0.\]
In particular, $P_{_F}(G_0)\leq P_{_F}(\Omega)$.	
\end{corollary}

\begin{proof}
Since $E_t \subset  G_t$ and, by Proposition \ref{min-hull} (b), $G_t$ is strictly outward $F$-minimizing, $E_t'\subset G_t$. Similarly, $\Omega'\subset G_0$.

We now argue by contradiction and assume that $|G_t \setminus E_t'|>0$ for $t>0$. Then, since $E_t'$ is strictly outward $F$-minimizing and $E_t'\subset G_t$, $P_{_F}(E_t') < P_{_F}(G_t)$. However, since $G_t\setminus E_t'\subset \{v=t\}$, Lemma \ref{v=t} (see also Corollary \ref{div=0}) implies $|Dv|_{_F}(G_t\setminus E_t')=0$ thus, by Proposition \ref{equiv-sol},
$$P_{_F}(G_t)\leq P_{_F}(E_t')+|Dv|_{_F}(G_t\setminus E_t')=P_{_F}(E_t'),$$
which contradicts $P_{_F}(E_t') < P_{_F}(G_t)$. Hence the statement follows.

Analogously, if we assume that $|G_0\setminus \Omega'|>0$, arguing as before (using also Remark \ref{Dv=0-on-v=t}) we reach a contradiction. 
\end{proof}

\section{Explicit examples of solutions}

{We begin with showing that Wulff shapes play the role of spheres in the case of the A-IMCF \eqref{weak-IMCF}; i.e. the level sets of the solution are Wulff shapes with exponentially growing radius.

\begin{example}\label{ex:1} \normalfont
  Let  $\Omega=\mathcal W_R$, for some $R>0$. Then, the solution to \eqref{EquationVp1} is given by $$v(x)=(N-1)\log\left(\frac{F^\circ(x)}{R}\right).$$
 Note that $v\in {\rm Lip}_{\rm loc}(\Omega^e)$, $\nabla v(x)=\frac{(N-1)\nabla F^\circ(x)}{F^\circ(x)}$ a.e. in $\Omega^e$ and $v(x)=0$ if $x\in\partial\Omega$. Then, by Lemma \ref{relationsF} (i), we have $F(\nabla v(x))=\frac{N-1}{F^\circ(x)}$.

 On the other hand, considering ${\mathbf z}(x):=\frac{x}{F^\circ(x)}$, it is immediate to check that ${\rm div}(\mathbf z)=\frac{N-1}{F^\circ(x)}$. Finally, since by Lemma \ref{relationsF}, ${\mathbf z}\in \partial F(\nabla F^\circ(x))$, we conclude that $v$ is the solution to \eqref{EquationVp1}.\end{example}

\begin{example} \normalfont
  Take a rectangle $\Omega=]-\frac{1}{2},\frac{1}{2}[\times ]-1,1[$ in $\R^2$ with norm $F = \ell_\infty$, and hence $\mathcal W_1=]-1,1[^2$. The solution is given by
   $$v(x,y)=\left\{\begin{array}
    {cc} \displaystyle \frac{1}{2}\log\left(\frac{|y|+\sqrt{y^2-3/4}}{3(|y|-\sqrt{y^2-3/4})}\right), & {\rm if \ } y^2>x^2+3/4, \medskip \\
     \displaystyle \frac{1}{2}\log\left(\frac{|x|+\sqrt{x^2+3/4}}{3(\sqrt{x^2+3/4}-|x|)}\right), & {\rm if \ } y^2<x^2+3/4.
  \end{array}\right.$$
  It suffices to define $$\z(x,y):=\left\{\begin{array}
    {cc} \displaystyle \left(\frac{x}{\sqrt{y^2-3/4}},{\rm sign}(y)\right), & {\rm if \ } y^2>x^2+3/4, \\
     \displaystyle \left({\rm sign }(x),\frac{y}{\sqrt{x^2+3/4}}\right), & {\rm if \ } y^2<x^2+3/4.
  \end{array}\right.$$
 Notice that the vertices of the level sets do not evolve linearly (see Figure \ref{fig-rectangle}), they instead move along the curve $y^2 = x^2+3/4$.

 \vspace*{-0.5cm}
    \begin{figure}[H]
    $$\includegraphics[scale=0.6]{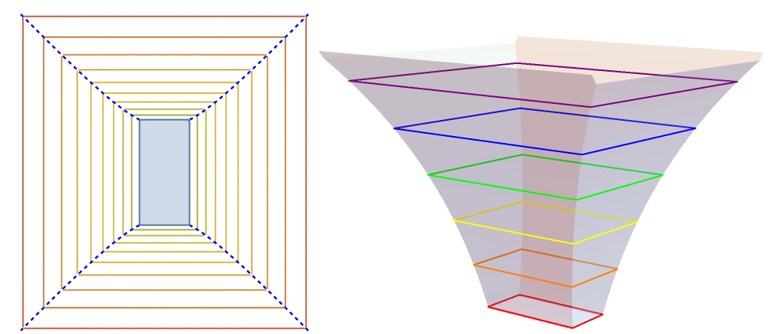}$$
\vspace*{-0.5cm}    \caption{Evolution of an initial rectangle with respect to the $\ell_1$-anisotropy. Plot of the level sets.}\label{fig-rectangle}
  \end{figure}
\end{example}

We finish this section with an example of fattening:

\begin{example}  \normalfont
  Let $\Omega\subset\R^2$ be the union of 3 squares centred at $(0,2),(\pm 2,-2)$ with sides of length $1$ and consider $F(x,y)=\|(x,y)\|_1$. In this case, $\Omega=\mathcal W_{\frac{1}{2}}(0,2)\cup \mathcal W_{\frac{1}{2}}(-2,-2)\cup \mathcal W_{\frac{1}{2}}(2,-2)$. We claim that the solution to \eqref{EquationVp1} is the following function (see Figure \ref{fig:3squares}):
  $$v(x,y)=\left\{\begin{array}{ll} \log(2), & {\rm if \ } (x,y)\in \mathcal W_{3}\setminus \overline\Omega, \\ \log (\|x,y-2\|_\infty), & {\rm if \ } (x,y)\in\mathcal W_{1}(0,2)\setminus \overline{\mathcal W_{\frac{1}{2}}(0,2)},\\ \log (2\|x+2,y+2\|_\infty), & {\rm if \ } (x,y)\in\mathcal W_{1}(-2,-2)\setminus \overline{\mathcal W_{\frac{1}{2}}(-2,2)}, \\ \log (2\|x-2,y+2\|_\infty), & {\rm if \ } (x,y)\in\mathcal W_{1}(2,-2)\setminus \overline{\mathcal W_{\frac{1}{2}}(2,2)}, \\ \log\left(\frac{2}{3}\|x,y\|_\infty\right), & {\rm if \ } (x,y)\in\R^2\setminus \overline{\mathcal W_{3}}.
  \end{array}\right.$$
  Notice that there is fattening at time $t^*=\log 2$: $E_{t^*}=\mathcal W_{1}(0,2)\cup\mathcal W_{1}(-2,2)\cup\mathcal W_{1}(2,-2)$ while $G_{t^*}=\mathcal W_3$. Moreover, it holds $P_F(G_{t^*})=P_F(E_{t^*})=24=e^{t^*}P_F(\Omega)$, consistent with \eqref{exp-per}.

  \vspace*{-0.6cm}
  \begin{figure}[H]
  	$$\includegraphics[width=0.8\textwidth]{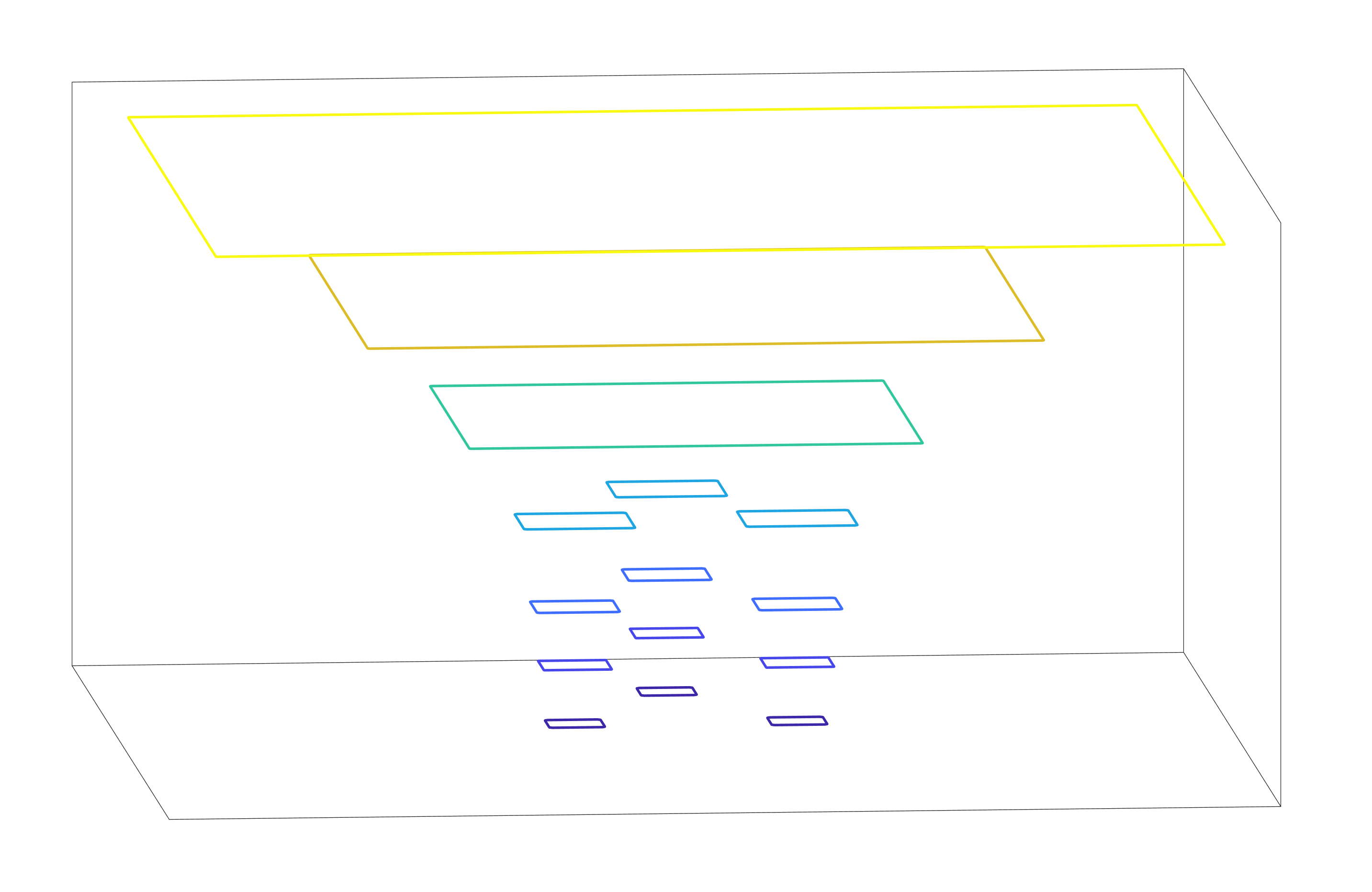}$$\vspace*{-1cm}\caption{Example of fattening.}\label{fig:3squares}
  \end{figure}

\vspace*{-0.4cm}

  A direct computation shows that $$\tilde \z(x,y):=\left\{\begin{array}{ll} \displaystyle \frac{(x,y-2)}{\|x,y-2\|_\infty}, & (x,y)\in\mathcal W_{1}(0,2)\setminus \overline{\mathcal W_{\frac{1}{2}}(0,2)}, \medskip \\ \displaystyle \frac{(x+2,y+2)}{\|x+2,y+2\|_\infty}, & (x,y)\in\mathcal W_{1}(-2,2)\setminus \overline{\mathcal W_{\frac{1}{2}}(-2,2)}, \medskip \\ \displaystyle \frac{(x-2,y+2)}{\|x-2,y+2\|_\infty}, & (x,y)\in\mathcal W_{1}(2,-2)\setminus \overline{\mathcal W_{\frac{1}{2}}(2,-2)}, \medskip \\ \displaystyle \frac{(x,y)}{\|x,y\|_\infty}, & (x,y)\in\R^2\setminus \overline {\mathcal W_{3}}.
  \end{array}\right.$$
  satisfies all the requirements for $v$ to be a solution in $(\R^2\setminus \overline{\mathcal W_3})\cup \left(\mathcal W_1(0,2)\cup\right.$ $\left.\mathcal W_1(-2,2)\cup \mathcal W_1(2,-2)\right)\setminus \overline\Omega$. In order to complete the proof, we need to extend $\tilde \z$ to $A:=\mathcal W_3\setminus\overline{(\mathcal W_1(0,2)\cup\mathcal W_1(-2,2)\cup \mathcal W_1(2,-2))}$. To achieve this aim, we need a vector field $\z^A: A\to \mathcal W_1$ such that
  $$\left\{\begin{array}{ll} {\rm div \ } \z^A=0 & {\rm in  \ } A, \smallskip \\ {\rm compatibility \ conditions } & {\rm on \ }\partial A. \end{array}\right. $$

  Next, we briefly sketch how to specify the latter on $\partial A$. Let $$\z:=\left\{\begin{array}{ll} \tilde \z & {\rm in \ } \R^2\setminus(\overline \Omega\cup A), \smallskip \\ \z^A=(z^A_1,z^A_2) & {\rm in \ }A. \end{array}\right.$$
  Then, since ${\rm div \,  \z}=F(Dv)<<\mathcal L^2(\R^2\setminus\overline\Omega)$, it follows that $\z^A\cdot\nu=\tilde\z\cdot \nu$, with $\nu$ being the outer unit normal to $\partial A$. Therefore, we need \begin{equation}\label{compatibility}\left\{\begin{array}{ll} z^A_1(\pm 3, y)=\pm 1, & z^A_1(\pm 1,y)=\pm {\rm sign}(y) \medskip \\ z^A_2(x,\pm 3)=1, & z^A_2(x,\pm 1)=\mp 1, \end{array}\right.\end{equation} for $(x,y)\in\partial A$. In Figure \ref{fig:arrows}, we plot the compatibility conditions at $\partial A$. 
  The arrows represent the orientation of the normal component of $\z^A$.

  The construction of  $\z^A$ follows as in \cite[Theorem 5]{BelletiniCasellesNovaga} (see \cite[Theorem 5]{LMM} for the $\ell_1$ anisotropy). In particular, it is obtained as a minimizer of
  $$\int_{A}({\rm div} \, {\boldsymbol\xi})^2\ d\mathcal L^2$$ among vector fields $\boldsymbol{\eta}\in L^\infty(A)$ satisfying  \eqref{compatibility}. Moreover, ${\rm div} \boldsymbol{\xi}={\rm div} \boldsymbol{\eta}$ is constant for any two minimizers $\boldsymbol{\xi}$, $\boldsymbol{\eta}$.
 Now, using integration by parts, $$\int_{A} {\rm div} \, \z^A \,d\mathcal L^2=\int_{\partial A}\z^A\cdot \nu \, ds=0,$$ and, since ${\rm div} \,\z^A$ is constant, necessarily ${\rm div} \, \z^A=0$. This finally shows that $v$ is the solution to \eqref{EquationVp1} with associated vector field $\z$.

  \vspace*{-0.5cm}
    \begin{figure}[H]
     $$\includegraphics[width=0.5\textwidth, scale=0.5]{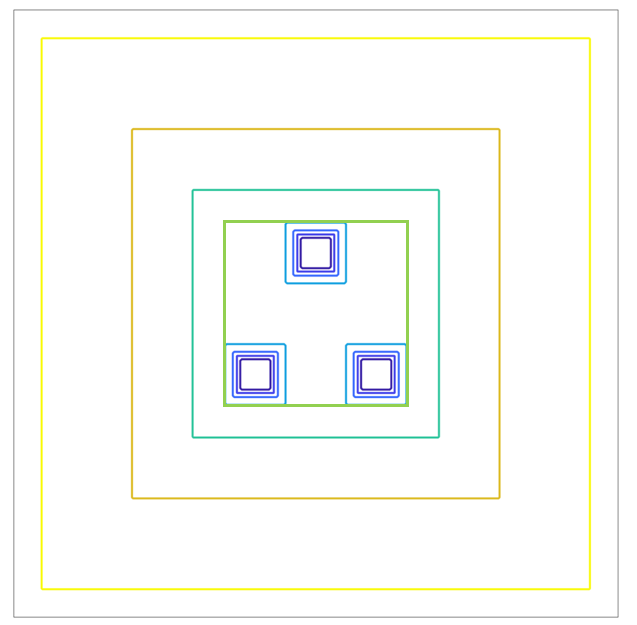}\includegraphics[width=0.5\textwidth]{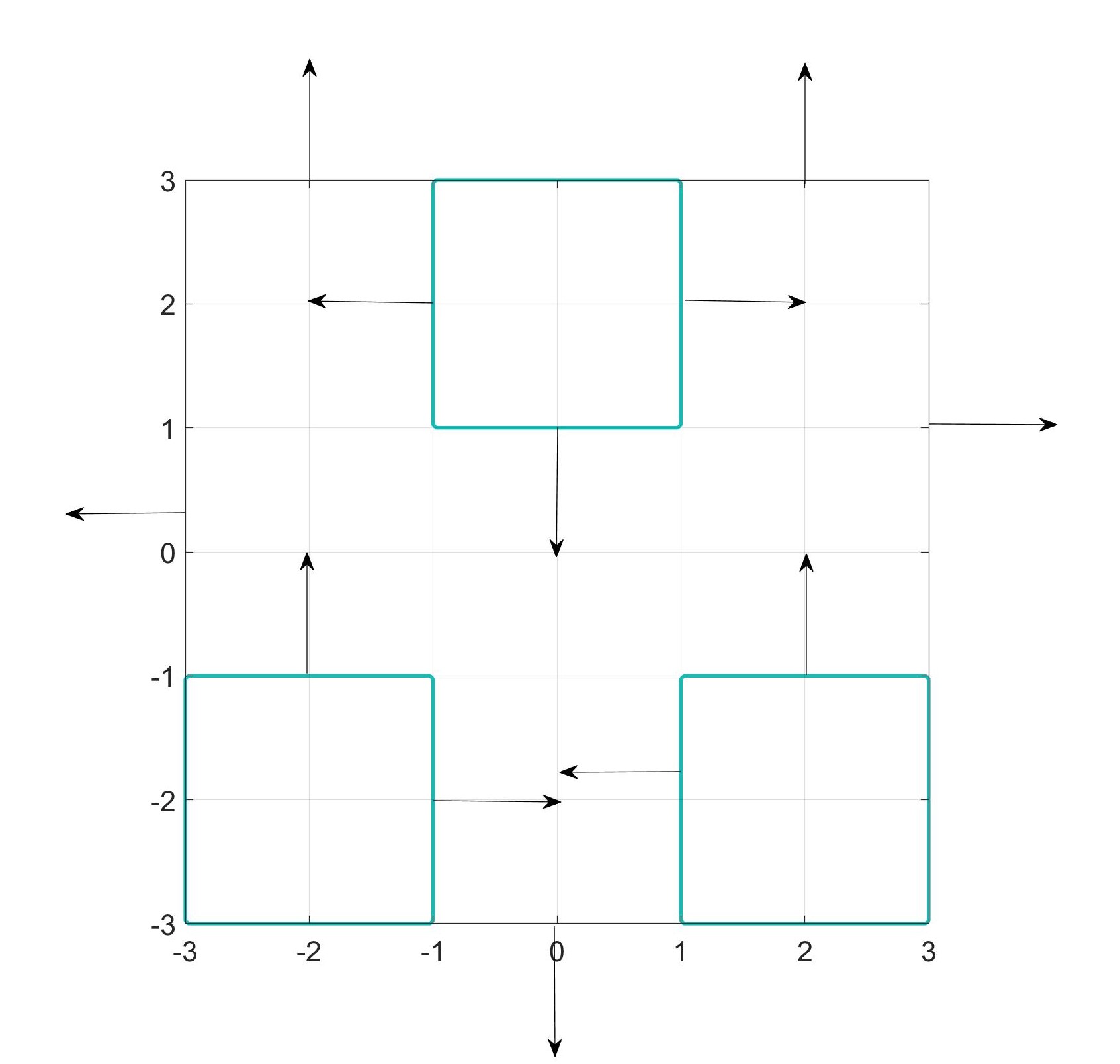}$$\vspace*{-0.5cm}\caption{Left: level sets of the solution. Right: Compatibility conditions of the vector field at the time of fattening.}\label{fig:arrows}
  \end{figure}
\end{example}


\begin{thebibliography}{00}
	
\bibitem{AmBe} M. Amar and G.\ Bellettini, A notion of total variation depending on a metric with discontinuous coefficients. {\it Ann. Inst. H. Poincaré C Anal. Non Linéaire} {\bf 11} (1994), no. 1, 91--133.

\bibitem{ACM} L.\ Ambrosio, G.\ Crippa, S.\ Maniglia, Traces and fine properties of a BD class of vector fields and applications. {\it Ann.\ Fac.\ Sci.\ Toulouse Math.\hspace*{-0.2cm} } (6) {\bf 14} (2005), no. 4, 527--561.	

	
\bibitem{Ambrosio} L. Ambrosio, N. Fusco, and D. Pallara, {\it Functions of bounded variation and free discontinuity problems}. Oxford Mathematical Monographs. The Clarendon Press, Oxford University Press, New York, 2000.





\bibitem{ABCM} F. Andreu, C. Ballester, V. Caselles and J.M. Mazón,  Minimizing total variation flow. {\it Differential Integral Equations} {\bf 14} (2001), no. 3, 321–360.

\bibitem{Andreu-Caselles-Mazon-ibero}
F. Andreu, V. Caselles and J. M. Mazón,
A Parabolic Quasilinear Problem for
Linear Growth Functionals.
{\it Rev. Mat. Iberoamericana}, {\bf 18} (2002), 135--185.

\bibitem{Anz} G.\ Anzellotti,
Pairings between measures and bounded functions and compensated compactness.
{\it Ann. Mat. Pura Appl.} (4) {\bf 135} (1983), 293–318 (1984).



\bibitem{Tama} R.C. Bassanezi and I. Tamanini,  Subsolutions to the least area problem and the "minimal hull" of a bounded set in $\R^n$. {\it Ann. Univ. Ferrara} {\bf 30} (1984), 27–40.

	\bibitem{BC} H. H. Bauschke, P. L. Combettes, {\it Convex Analysis
	and Monotone Operator Theory in Hilbert Spaces}. Second Edition. CMS Books in Mathematics, Springer, 2011.

\bibitem{BelletiniCasellesNovaga} G. Bellettini, V. Caselles and M- Novaga,  The total variation flow in $\R^N$.
    {\it J. Differential Equations} (2) {\bf 184} (2002), 475-525.



\bibitem{BHW} S. Brendle, P.-K. Hung, M.-T. Wang, A Minkowski inequality for hypersurfaces in the Anti-deSitter–Schwarzschild manifold, {\it Comm. \hspace*{-0.2cm} Pure Appl.\hspace*{-0.05cm} Math.} \hspace*{-0.05cm} {\bf 69}(1) (2016) 124--144.

\bibitem{BW} S. Brendle, M.-T. Wang, A Gibbons–Penrose inequality for surfaces in Schwarzschild spacetime,
{\it Comm. Math. Phys.} {\bf 330} (1) (2014) 33--43.
	
\bibitem{CRMS} E. Cabezas-Rivas, S. Moll and M. Solera, Characterization of the subdifferential and minimizers for the anisotropic p-capacity, {\it to appear in Adv. Calc. Var.}, \url{arxiv:2305.03498} (2023).



\bibitem{CFM} V.\ Caselles, G.\ Facciolo and E.\ Meinhardt,
Anisotropic Cheeger sets and applications. {\it SIAM J. Imaging Sci.} {\bf 2} (2009), no. 4, 1211–1254.

\bibitem{Cas} V. Caselles, On the entropy conditions for some flux limited diffusion equations. J. Differential Equations {\bf 250} (2011), no. 8, 3311--3348.

\bibitem{Cha-Nov} A. Chambolle and M. Novaga, Anisotropic and crystalline mean curvature flow of mean-convex sets, {\it Ann. Sc. Norm. Super. Pisa Cl. Sci. (5)} {\bf 23} (2022), no. 2, 623–643.

\bibitem{CoCraCiMa} G.E. Comi, G. Crasta, V. De Cicco and A Malusa, Representation formulas for pairings between divergence-measure fields and  functions, {\it J. Funct. Anal.} {\bf 286} (2024), no.1, paper no. 110192, 32 pp.


\bibitem{CraCi} G.\ Crasta and V.\ De Cicco, Anzellotti's pairing theory and the Gauss-Green theorem. {\it Adv.\ Math.} {\bf 343} (2019), 935--970.

\bibitem{CheFri} G-Q.\ Chen and H.\ Frid, Divergence-measure fields and hyperbolic conservation laws. {\it Arch. Ration. Mech. Anal.} {\bf 147} (1999), no. 2, 89--118.

\bibitem{IAMCF} F. Della Pietra, N. Gavitone and C. Xia, {\it Motion of level sets by inverse anisotropic mean curvature}, {\it Commun. Anal. Geom.} {\bf 31}, (2023), no. 1, 97–118.
	
	\bibitem{Evans} L.C. Evans, {\it Partial Differential Equations: Second Edition}. American Mathematical Soc., 2010.



	
\bibitem{Gerhardt} C. Gerhardt,
	 Flow of nonconvex hypersurfaces into spheres.
	{\it J. Differential Geom.} {\bf 32} (1990), no. 1, 299--314.
	
\bibitem{GMP} L. Giacomelli, S. Moll, F. Petitta, Nonlinear diffusion in transparent media: the resolvent equation, {\it Adv. Calc. Var.} {\bf 11} (4) (2018), 405–432.
	

	
\bibitem{HuIl} G. Huisken and T. Ilmanen,
{\it The inverse mean curvature flow and the Riemannian Penrose inequality.}
J. Differential Geom. {\bf 59} (2001), no. 3, 353–-437.

\bibitem{HuIl2}  G. Huisken, T. Ilmanen, Higher regularity of the inverse mean curvature flow, {\it J. Differential Geom.} {\bf 80} (3) (2008) 433–451.

\bibitem{JN} F. John and L. Nirenberg, {\it On functions of bounded mean oscillation.} Comm. Pure Appl. Math. {\bf 14} (1961), 415--426.

\bibitem{kohn-temam}
R. Kohn and R. Temam,
{Dual space of stress and strains with applications to Hencky plasticity,}
{\it Appl. Math. Optim.}, {\bf 10} (1983) 1--35.

\bibitem{KoNi} B. Kotschwar and L. Ni. Local gradient estimates of $p$-harmonic functions, $1/H$-flow, and an entropy formula. {\it Ann Sci Éc Norm Supér (4)}, {\bf  42} (2009) no. 1, 1--36.

\bibitem{LeSa} Leonardi, G.P., Saracco, G. The prescribed mean curvature equation in weakly regular domains. {\it Nonlinear Differ. Equ. Appl.} {\bf 25}, 9 (2018).

    	\bibitem{Leoni} G. Leoni, {\it A first course in Sobolev spaces}. Second edition. Graduate Studies in Mathematics, 181. American Mathematical Society, Providence, RI, 2017.

\bibitem{LMM} M. Łasica, S. Moll and P. Mucha,
Total variation denoising in l1 anisotropy. {\it SIAM J. Imaging Sci.} {\bf 10} (2017), no.4, 1691–1723.




\bibitem{Maggi} F. Maggi,  {\it Sets of finite perimeter and geometric variational problems (an introduction to geometric measure theory)}. Cambridge Studies in Advanced Mathematics (2012)


\bibitem{MazSer13} J.M. Mazón and S. Segura de León, The Dirichlet problem for a singular elliptic equation arising in the level set formulation of the inverse mean curvature flow. {\it Adv. Calc. Var.} {\bf 6} (2013), no. 2, 123--164.


\bibitem{MazSer15} J.M. Mazón and S. Segura de León, A non-homogeneous elliptic problem dealing with the level set formulation of the inverse mean curvature flow. {\it J. Differential Equations} {\bf 259} (2015), no. 7, 2762--2806.

\bibitem{moll_05}
J. S. Moll, The anisotropic total variation flow.
{\it Math. Ann.} 332,(2005),  177–218 .

	\bibitem{Moser} R. Moser, The inverse mean curvature flow and p-harmonic functions. {\it J. Eur. Math. Soc. (JEMS)} {\bf 9} (2007), no. 1, 77--83.
	
	\bibitem{Moser08} R.\ Moser,  The inverse mean curvature flow as an obstacle problem. {\it Indiana Univ. Math. J.} {\bf 57} (2008), no. 5, 2235–2256.
	
	\bibitem{Moser15} R.\ Moser, Geroch monotonicity and the construction of weak solutions of the inverse mean curvature flow. {\it Asian J. Math.} {\bf 19} (2015), no. 2, 357–376.



\bibitem{Smo} K. Smoczyk, Remarks on the inverse mean curvature flow, {\it Asian J. Math.} {\bf 4}(2) (2000), 331--335.

\bibitem{Stam} G. Stampacchia, Équations elliptiques du second ordre à coefficients discontinus, {\it Séminaire Jean Leray}, no. 3, 1--77 (1963–1964)

\bibitem{Stu} F. Schuricht, A new mathematical foundation for contact interactions in continuum physics, {\it Arch. Ration. Mech. Anal.} {\bf 184} (3) (2007) 495–551.

\bibitem{Urbas} J. Urbas, On the Expansion of Star-shaped Hypersurfaces by Symmetric Functions of Their Principal Curvatures, {\it Math. Z.} {\bf 205} (1990), 355--372.

\bibitem{Xia} C. Xia, {\it On a class of anisotropic problem}. PhD Thesis,  University Freiburg,
	2012.
	
\bibitem{Xia2} C. Xia,
Inverse anisotropic mean curvature flow and a Minkowski type inequality,
{\it Advances in Mathematics}, {\bf 315} (2017), 102--129.
	

	
	
	

	
	
	
	
\end{thebibliography}
\end{document}